\documentclass[11pt]{article}

\usepackage{color}
\usepackage{mathptmx}       
\usepackage{helvet}         
\usepackage{courier}        
\usepackage{type1cm}        
\usepackage{graphicx}        
\usepackage[bottom]{footmisc}

\usepackage{amsmath}    
\usepackage{amssymb}    
\usepackage{amsthm}    
\usepackage{epsfig}    
\usepackage{psfrag}     
\usepackage{mathrsfs}  

\newtheorem{theorem}{Theorem}[section]

\newtheorem{definition}[theorem]{Definition}

\newtheorem{proposition}[theorem]{Proposition}

\newtheorem{remark}[theorem]{Remark}

\newcommand\DT[1]{\mathchoice
                 {{\buildrel{\hspace*{.1em}\text{\LARGE.}}\over{#1}}}
                 {{\buildrel{\hspace*{.1em}\text{\Large.}}\over{#1}}}
                 {{\buildrel{\hspace*{.1em}\text{\large.}}\over{#1}}}
                 {{\buildrel{\hspace*{.1em}\text{\large.}}\over{#1}}}}
\newcommand\DDT[1]{\mathchoice
   {{\buildrel{\hspace*{.13em}\text{\LARGE.\hspace*{-.13em}.}}\over{#1}}}
   {{\buildrel{\hspace*{.1em}\text{\Large.\hspace*{-.1em}.}}\over{#1}}}
   {{\buildrel{\hspace*{.1em}\text{\large.\hspace*{-.1em}.}}\over{#1}}}
   {{\buildrel{\hspace*{.1em}\text{\large.\hspace*{-.1em}.}}\over{#1}}}}
\newcommand{\lineunder}[2]{\LU{\begin{array}[t]{c}\underbrace{#1}\vspace*{.5em}\end{array}}{\mbox{\footnotesize\rm #2}}}
\newcommand{\linesunder}[3]{\LSU{\begin{array}[t]{c}\underbrace{#1}\vspace*{.5em}\end{array}}{\mbox{\footnotesize\rm #2}}{\mbox{\footnotesize\rm#3}}}
\newcommand{\LU}[2]{\begin{array}[t]{c}#1\vspace*{-1em}\\_{#2}\end{array}}
\newcommand{\LSU}[3]{\begin{array}[t]{c}#1\vspace*{-1em}\\_{#2}\vspace*{-.5em}\\_{#3}\end{array}}
\newcommand{\SEFS}[1]{\sigma_\text{\tiny{DAM}}^{#1}}

\newcommand{\overlinealpha}{\hspace{.1em}\overline{\hspace{-.1em}\alpha}}
\newcommand{\underlinealpha}{\hspace{.1em}\underline{\hspace{-.1em}\alpha}}
\newcommand{\COL}[1]{{\color{black}#1}}
\newcommand{\NEW}[1]{{\color{black}#1}}

\allowdisplaybreaks
\newcommand{\wt}{\widetilde}

\newcommand{\COMMENT}[1]{{\color{red}\bf}}
\newcommand{\QUE}[1]{{\color{red}\bf}}

\begin{document}
\begin{sloppypar}

\begin{center}
    \vspace*{-1em}{\small\tt
      A preprint to appear in:
    ''Topics in Applied Analysis and Optimisation''
    \\[-.0em](Eds. J.-F. Rodrigues and M. Hintterm\"uller)
    CIM Series in Math. Sci., Springer.
    \\[.3em]\hrule\hrule\hrule} 
    \vspace*{5em}
  {\LARGE\bf
Models of dynamic damage
and phase-field fracture,\\[.4em]and their various time discretisations}


  \vspace*{3em}

    {\sc Tom\'{a}\v{s} Roub\'\i\v{c}ek}\\[1em]
Charles University,
Mathematical Institute\\ Sokolovsk\'a 83, CZ-186~75~Praha~8,
Czech Republic\\
and\\
Institute of Thermomechanics, Czech Academy of Sciences\\Dolej\v skova~5,
CZ-182~08 Praha 8, Czech Republic\\
tomas.roubicek@mff.cuni.cz





\end{center}

\def\fDAM{{\phi}}
\def\fDAMprime{{\phi'}}
\def\Frakg{{\mathfrak{g}}}
\def\eps{\varepsilon}
\def\eq{\eqref}
\def\R{{\mathbb R}}
\def\N{{\mathbb N}}
\def\Rsym{\R^{d\times d}_{\rm sym}}
\def\bbD{{\mathbb D}}
\def\bbC{{\mathbb C}}
\def\bbI{{\mathbb I}}
\def\Colon{{\colon}}
\def\Cdot{\!\cdot\!}
\def\d{{\rm d}}
\def\pl{{\partial}}
\def\taur{\chi}
\def\In{\!\in\!}
\def\zetai{\zeta_\text{\sc i}}
\def\zetav{\zeta_\text{\sc v}}
\def\SHEAR{G}
\def\SYLD{\sigma_\text{\sc yld}}

\abstract{Several variants of models of damage in viscoelastic continua
  under small strains in the Kelvin-Voigt rheology \COL{are} presented and
  analyzed
by using the Galerkin method. The particular case, known as a phase-field 
fracture approximation of cracks, is discussed in detail. All these models 
are dynamic (i.e.\ involve inertia to model vibrations or waves possibly 
emitted during fast damage/fracture or induced by fast varying forcing) and 
consider viscosity which is also damageable. Then various options for time 
discretisation are devised. Eventually, extensions to more complex rheologies 
or a modification for large strains are briefly exposed, too.}

\section{Introduction}
Damage in continuum mechanics of solids is an important part of engineering 
modelling (and also experimental research), focusing on the attribute of 
various degradation of materials. During past few decades, some of 
the engineering models had been also under rigorous mathematical scrutiny.

Phenomenological damage models structurally represent the simplest example of 
the concept of internal variables, where only one 
scalar-valued variable (here
denoted by $\alpha$) is considered. Cf.\ G.A.\,Maugin \cite{Maug15SIVS} 
for a thorough historical survey of this concept.
This scalar-phenomenological-damage concept was invented by 
L.M.\,Kachanov \cite{Kach58TRPD} \COL{and Yu.N.\,Rabotnov \cite{Rabo69CPSM}},
the damage variable ranging the interval $[0,1]$ and having an intuitive
microscopical interpretation as a density of microcracks or microvoids.
There are two conventions: damaging means $\alpha$ increasing and $\alpha=1$
means maximal damage (which is used in engineering
or e.g.\ also in geophysics) or, conversely, damaging means $\alpha$ decreasing
and $\alpha=0$ means maximal damage (which is used in mathematical
literature and also here)\NEW{, cf.\ e.g.\ the monographs
  \cite[Ch.12]{Frem02NST} and \cite[Ch.6]{Frem12PCM}}. 
{Let us still note that, although 
damage as a single variable is most often used in applications, some 
models with more variables are sometimes considered in engineering, too. 
Most generally, one may think about 8th-order tensor as a damage variable, 
transforming 4th-order elastic-moduli tensor $\bbC$, cf.\ e.g.\ 
\cite{Mura12CDM}.}

Damage can be (and typically is) a very fast process, usually much faster
than the time scale of external loading. This is reflected by an (often
accepted) idealization to model it as a rate-independent process which
can have arbitrary speed. No matter whether the model is rate-independent
or involves some sort of damage viscosity, the fast damage may generate 
elastic waves in the continuum. Conversely, waves can trigger damage.
This combination of damage at usually localized areas and inertial effects 
in the whole bulk needs a bit special methods both for
rigorous analysis and for numerical approximation, some of them being suitable
rather for vibration (where transfer of kinetic energy is not dominant) 
than waves. On top of it, in some application\COL{s} even the loading itself
can \COL{vary} fast in time during various impacts or explosions.
Such {\it dynamic damage}\index{dynamic damage/fracture} or  dynamic
{\it fracture mechanics} \cite{Freu98DFM} is
the main focus of this chapter, the (often considered) quasistatic variants 
being thus intentionally avoided here. It is also important that inertia 
suppresses artificial global long-range interactions which otherwise make 
various unphysical effects and causes a need of some rather artificial 
quasistatic models, cf.~Remark\ref{rem-FFM} below. \COMMENT{??????}

Most of this exposition will be formulated at small strains. 
Plain damage will be presented in several variants in 
Section~\ref{sect-damage-small-strains}. Its usage for fracture mechanics
exploiting the so-called phase-field approximation will then be in 
Section~\ref{sect-damage-PF-fracture}, outlining a wide menagerie 
of models towards distinguishing crack initiation and crack propagation,
possibly sensitive to modes (i.e.\ opening versus shearing), some of these
models being likely new.
Various discretisation\COL{s} of these models in time may exhibit
various useful properties, which will be presented in 
Section~\ref{sect-damage-discretisation}. Eventually, the basic 
scenario of small-strain models with just one scalar-valued damage variable 
can be enriched in many ways, by involving some other internal variables
like plastic strain or a \COL{diffusant} content and also temperature, which is 
certainly motivated by specific applications. One can make it
either in the framework of small strains considered in the previous
sections, or even at large strains. Some of these enhancements  
will be briefly outlined in Section~\ref{sect-damage-generalization}

\section{Models of damage at small strains}\label{sect-damage-small-strains}

Beside \COL{the} already mentioned alternative of damage 
being rate-dependent versus rate-independent, there are many variants.
Basic alternatives are unidirectional damage (i.e.\ no healing is allowed, 
relevant in most engineering materials) versus reversible damage (i.e.\ a 
certain reconstruction of the material is possible, relevant e.g.\ in 
rock mechanics in the time scales of thousands years or more). 
And, of course, damage models can be incorporated into various viscoelastic
models, and damage can influence not only the stored energy but also
the dissipation potential. The damage can be complete (which is 
mathematically much more difficult, cf.\ \cite{MiRoZe10CDEV}
at least for some partial results) or incomplete.
In addition to the simplest mode-insensitive damage, many applications
need a mode-sensitive damage (damage by tension/opening
easier than by shearing).
On top of all this, 
there is a conceptual discussion whether rather stress or energy 
(or a combination of both) causes damage.


In addition to these options, some nonlocal theories are typically used. 
This concern the damage variable and sometimes also the strain. 
Here we have in mind so-called weakly nonlocal concepts which involve
usual local gradients. The former case thus involves damage gradient into
the stored energy and allows us to introduce length-scale into the damage,
while the latter option allows for weaker assumptions on lower-order terms
and for involving dispersion into elastic waves, as \COL{discussed} in 
\cite{Jira04NTCM}.

There are many options of damage models outlined above, some of them  
complying with rigorous analysis while some others which making troubles.
Most mathematical models at small strains consider the specific stored energy 
$\varphi=\varphi(e,\alpha)$ quadratic in terms of the small-strain
variable $e\in\Rsym$. 

From an abstract viewpoint, the evolution is governed by {\it Hamilton's 
variational principle}\index{Hamilton variational principle} generalized for 
dissipative systems \cite{Bedf85HPCM}, which says that, among all admissible 
motions $q=q(t)$ on a fixed time interval $[0,T]$, the actual motion 
makes 
\begin{align}\label{hamilton}
\int_0^T{\mathscr L}\big(t,q,\DT{q}\big)\,\d t\ \ \mbox{ stationary\ \ 
(i.e.~$q$ is its 
\emph{critical point})},\index{critical point!for Lagrangian}
\end{align}
where $\DT{q}=\frac{\pl}{\pl t}{q}$ and ${\mathscr L}(t,q,\DT{q})$ is the 
\emph{Lagrangian}\index{Lagrangian} defined by 
\begin{align}\label{Lagrangean}
&{\mathscr L}\big(t,q,\DT q\big):={\mathscr T}\big(\DT q\big)
-{\mathscr E}(t,q)+\langle
F(t),q\rangle\,,
\end{align}
where $F=-\pl_{\DT q}{\mathscr R}(q,\DT q)$ is a nonconservative 
force assumed for a moment fixed, with ${\mathscr R}(q,\cdot)$ denoting the 
(Rayleigh's \emph{pseudo})\emph{potential}\index{potential!of dissipative force}
 of the {\it dissipative force}.
Then (\ref{hamilton}) leads after by-part integration in time to 
\begin{align}
\pl_q^{}{\mathscr L}\big(t,q,\DT{q}\big)
-\frac{\d}{\d t}\pl_{\DT q}^{}{\mathscr L}\big(t,q,\DT{q}\big)=0.
\end{align}
This gives the 
abstract 2nd-order evolution equation 
\begin{align}\label{abstract-momentum-eq}
{\mathscr T}'\DDT q+\pl_{\DT q}{\mathscr R}(q,\DT q)+\pl_q{\mathscr E}(t,q)=0
\end{align} 
where the apostrophe (or $\pl$) indicates the (partial) G\^ateaux differential.

In the context of this section, the state $q=(u,\alpha)$ consists from 
the displacement $\Omega\to\R^d$ and the damage profile 
$\Omega\to[0,1]$ with $\Omega\subset\R^d$ a bounded domain with a 
Lipschitz boundary $\Gamma$, and we specify the overall
kinetic energy, stored energy (including external loading), and 
dissipation potential as
\begin{subequations}\label{damage-gradient}\begin{align}
&\label{kinetic-energy}
{\mathscr T}\big(\DT q\big)={\mathscr T}\big(\DT u\big):=
\int_\Omega\frac\varrho2\big|\DT u\big|^2\,\d x\,,
\\&
{\mathscr E}(t,q)={\mathscr E}(t,u,\alpha):=
\int_\Omega
\varphi(e(u),\alpha)+\frac\kappa{\COL{2}}|\nabla\alpha|^2
\label{damage-gradient-E}
{-}f(t)\cdot u
+\delta_{[0,1]}^{}(\alpha)\,\d x+\int_\Gamma g(t)\cdot u\,\d S\,,
\\&{\mathscr R}(q,\DT q)={\mathscr R}(\alpha,\DT u,\DT\alpha)
=\int_\Omega\frac12\bbD(\alpha)e(\COL{\DT u})\Colon e(\COL{\DT u})+\zeta(\DT\alpha)\,\d x
\end{align}\end{subequations}
with the small-strain tensor $e(u)=\frac12(\nabla u^\top\!{+}\nabla u)$ and with 
some specific damage dissipation-force potential 
$\zeta:\R\to[0,+\infty]$ convex with $\zeta(0)=0$, 
with a 4th order tensor $\bbD:[0,1]\to\R^{d\times d\times d\times d}$ 
smoothly dependent on $\alpha$, \NEW{$\kappa>0$ a phenomenological
  coefficient determining a length-scale of damage} \COL{(which is a usual
  engineering concept, cf.\ e.g.\
  \cite{BazJir02NIFP}, useful also from analytical reasons)}, and with 
$\delta_{[0,1]}^{}(\cdot):\R\to\{0,\NEW{+}\infty\}$ denoting
the indicator function of the interval $[0,1]$ where the damage variable is
assumed to take its values.

This general framework gives 
a relatively simple model of damage in the linear 
Kelvin-Voigt viscoelastic solid\COL{s} where the only internal variable 
is the damage. 
Feeding \eqref{abstract-momentum-eq} by the 
functionals \eqref{damage-gradient}, 
we arrive at the system of partial differential equation and inclusion\COL{s} 
\begin{subequations}\label{8-GSM1**}\begin{align}\label{8-GSM1**-u}
&\varrho\DDT{u}
-{\rm div}\big(\bbD(\alpha)e(\DT{u})+
\pl_e^{}\varphi(e(u),\alpha)
\big)
=f&&\text{in }Q,\\\label{8-GSM1**-z}
&\pl\zeta(\DT\alpha)+
\pl_\alpha^{}\varphi(e(u),\alpha)
-{\rm div}\big(\kappa
\nabla \alpha\big)
+r_{\text{\sc c}}^{}\ni0
&&\text{in }Q,
\\&\label{8-GSM1**-r}
r_{\text{\sc c}}^{}\in N_{[0,1]}(\alpha)&&\text{in }Q,
\intertext{where $N_{[0,1]}=\pl\delta_{[0,1]}$ is the normal cone to the interval 
$[0,1]$ where $\alpha$ is supposed to be valued, together with the
boundary conditions}
&\label{8-GSM1**BC}
\big(\bbD(\alpha)e(\DT{u})+
\pl_e^{}\varphi(e(u),\alpha)
\big)\vec{n}=g\ \ \text{ and }\ \ \
\nabla \alpha\Cdot\vec{n}=0&&\text{on }\Sigma,
\end{align}\end{subequations}
where $Q=\Omega\times I$ and $\Sigma=\Gamma\times I$ with $I=[0,T]$ for 
a fixed time horizon $T>0$, and where $\vec{n}$ is the outward unit 
normal to $\Gamma$. In fact, (\ref{8-GSM1**}b,c) can be understood as 
one doubly-nonlinear inclusion if the ``reaction pressure'' $r_{\text{\sc c}}^{}$ 
would be substituted from \eq{8-GSM1**-z} into \eq{8-GSM1**-r}.
We will consider an initial-value problem and thus 
complete \eqref{8-GSM1**} by the initial conditions 
\begin{align}\label{8-GSM1**IC}
u|_{t=0}^{}=u_0,\ \ \ \ \ \DT u|_{t=0}^{}=v_0,\ \ \ \ \ \alpha|_{t=0}^{}=\alpha_0\ \ \ \ \
\text{ \NEW{in} }\ \Omega. 
\end{align}
 The energetics can be obtained by testing \eq{8-GSM1**-u} by $\DT u$ and 
 \eq{8-GSM1**-z} by $\DT\alpha$. \COL{After integration over $\Omega$ with
   using Green's formula and by-part integration over a time interval $[0,t]$},
 this test yields, at least formally,\footnote{\COL{This means that
   \eq{KV-damage-engr} can rigorously be proved only for sufficiently smooth
   solutions, e.g.\ $\DT\alpha$ is to be in duality with ${\rm div}(\kappa\nabla\alpha)$, as e.g.\ in Proposition~\ref{prop-2+} below.}}
\begin{align}\nonumber
&\int_\Omega\!\!\!\!\!\linesunder{\frac\varrho2|\DT u(t)|^2}{kinetic}{energy at time $t$}\!\!\!\!\!+
\lineunder{\varphi(e(u(t)),\alpha(t))
+\frac\kappa 2|\nabla\alpha(t)|^2}{stored energy at time $t$}
\,\d x
+\int_0^t\!\!\int_\Omega\lineunder{\bbD(\alpha) e(\DT u)\Colon e(\DT u)
+\DT\alpha\pl\zeta(\DT\alpha)}{dissipation rate}\,\d x\d t
\\[-.4em]&
\qquad
=\int_\Omega\!\!\!\!\!\!\linesunder{\frac\varrho2|v_0|^2}{initial}{kinetic energy}\!\!\!\!\!\!+\lineunder{\varphi(e(u_0),\alpha_0)
+\frac\kappa 2|\nabla\alpha_0|^2}{initial stored energy}\,\d x
+\int_0^t\!\!\int_\Omega\!\!\!\linesunder{f\Cdot\DT u}{power of}{bulk load}\!\!\!\d x\d t+
\int_0^t\!\!\int_\Gamma\!\!\!\!\!\!\linesunder{g\Cdot\DT u}{power of}{surface load}\!\!\!\!\!\!\d x\d t.
\label{KV-damage-engr}
\end{align}

In fact, the model  \eq{8-GSM1**} may simplify in some particular situations
when  $r_{\text{\sc c}}^{}=0$ and \eqref{8-GSM1**-r} can be omitted, in
particular when
\begin{align}\label{damage-avoiding-[0,1]-constraints}
\pl_\alpha\varphi(e,0)\le0\ \ \text{ and }\ \ \
\begin{cases}
\pl_\alpha\varphi(e,1)\ge0
,\ \ \text{ or}\\
\zeta(\DT\alpha)=+\infty\ \ \text{ for }\ \DT\alpha>0\,.
\end{cases}
\end{align}
The first option allows for healing if $\zeta$ is finite also for $\DT\alpha>0$,
while the second option is called unidirectional damage. The condition 
$\pl_\alpha\varphi(e,0)=0$ needs infinitely large driving force to achieve
$\alpha=0$\NEW{, i.e.\ some sort of large hardening when damaging evolves}. 


The weak formulation of \eq{8-GSM1**-u} with the initial/boundary conditions
from \eq{8-GSM1**BC}-\eq{8-GSM1**IC} is quite standard, using usually 
one Green formula in space and one or two by-part integration\COL{s} in time.
The weak formulation of (\ref{8-GSM1**}b,c) consists in two variational 
inequalities. Writing the convex subdifferential in  \eq{8-GSM1**-z},
one see the term $\COL{\int_Q\DT{\alpha}\pl_\alpha\varphi(e(u),\alpha)}$
which is not a-priori integrable and we substitute \COL{it}
by using the calculus
\begin{align}
\int_Q\DT{\alpha}\pl_\alpha\varphi(e(u),\alpha)\,\d x\d t
=\int_\Omega\varphi\big(e(u(T)),\alpha(T)\big)-\varphi\big(e(u_0),\alpha_0\big)\,\d x
\label{chain-rule-for-GSM}
-\int_Q\pl_e\varphi(e(u),\alpha)\Cdot e(\DT{u})\,\d x\d t.
\end{align}
Thus, using the standard notation $L^p$, $W^{k,p}$, and $L^p(I;\cdot)$ 
or $W^{1,p}(I;\cdot)$ for Lebesgue, Sobolev, and Bochner or Bochner-Sobolev
spaces using also the convention $H^k:=W^{k,2}$, we arrive at:

\begin{definition}[Weak formulation]\label{def1}
A 
\COL{triple} $(u,\alpha,r_{\text{\sc c}}^{}){\in} 
H^1(I;H^1(\Omega;\R^d)){\times}H^1(Q){\times}$ ${\times}L^2(Q)$ is called a weak solution
to the initial-boundary-value problem \eq{8-GSM1**}--\eq{8-GSM1**IC} if 
$u|_{\COL{t=0}}^{}=u_0$ \NEW{and $0\le\alpha\le1$ hold a.e.}\ together with 
\begin{subequations}\begin{align}
&\hspace{0em}
\int_Q\!\big(\bbD(\alpha)e(\DT u)+\pl_e\varphi(e(u),\alpha)\big)\Colon e(v)
-\varrho \DT u\Cdot\DT{v}\,\d x\d t
=\!\int_\Omega\!v_0\Cdot v(0,\cdot)
\,\d x+\!\int_Qf\Cdot v\,\d x\d t
+\!\int_{\COL{\Sigma}}\!g\Cdot v\,\d S\d t
\label{7-very-weak-sln2}
\intertext{for all $v\In L^2(I;H^1(\Omega;\R^d))\,\cap\,\COL{H^{1}}(I;L^2(\Omega;\R^d))$
with $v|_{t=T}^{}=\DT{v}|_{t=T}^{}=0$, and}
\nonumber&\int_Q\!\pl_\alpha\varphi(e(u),\alpha)z
+r_{\text{\sc c}}^{}(z{-}\DT\alpha)
+\kappa
\nabla\alpha{\cdot}\nabla z+\zeta(z)\,\d x\d t
+\int_\Omega\varphi\big(e(u_0),\alpha_0\big)
+\frac\kappa 2|\nabla\COL{\alpha_0}|^2\,\d x
\\[-.4em]&\quad
\ge\!\int_Q
\zeta(\DT\alpha)+
\pl_e\varphi(e(u),\alpha)
\Colon e(\DT{u})\,\d x\d t
+\int_\Omega\varphi\big(e(u(T)),\alpha(T)\big)
+\frac\kappa 2|\nabla\COL{\alpha(T)}|^2\,\d x\,
\label{id6-}
\end{align}\end{subequations}
to be valid for all $z\in C^1(Q)$ and with $r_{\text{\sc c}}^{}$ 
satisfying \eq{8-GSM1**-r} a.e.\ on $Q$.
\end{definition}



Let us now analyze the model with the (partly) damageable viscosity 
in the special situation that
$\bbD(\alpha)=\bbD_0+\taur\pl_e\varphi(\cdot,\alpha)$ with a 
relaxation time $\taur>0$ possibly dependent on $x\In\Omega$, 
cf.\ \cite{LaOrSu10ESRM} or also \cite[Sect.5.1.1 and 5.2.5]{MieRou15RIST} 
for the rate-independent unidirectional damage. 
This means that $\varphi(\cdot,\alpha)$ is quadratic
and we thus specify the stored-energy density 
$\varphi:\Rsym\times[0,1]\to\R$ as
\begin{align}\label{damage-gradient-a}
&\varphi(e,\alpha)=
\frac12\bbC(\alpha)e\Colon e-\fDAM(\alpha)
\end{align}
with a 4th order elastic-moduli tensor $\bbC:[0,1]\to\R^{d\times d\times d\times d}$ 
continuously dependent on $\alpha$ and
with $\fDAM$ standing for the specific energy of damage which (if $\fDAM$ 
is increasing) gives rise to a driving force for healing.

This specifies the system (\ref{8-GSM1**}a-c) as
\begin{subequations}\label{8-GSM1*-special}
\begin{align}\label{8-GSM1**-u-special}
&\DT u=v,\ \ \ \ \ \ \ \ \ \varrho\DT{v}
-{\rm div}\big(\bbD(\alpha)e(v)+
\bbC(\alpha)e(u)
\big)
=f&&\text{in }Q,\\\label{8-GSM1**-z-special}
&\pl\zeta(\DT\alpha)+
\frac12\bbC'(\alpha)e(u)\Colon e(u)
-{\rm div}\big(\kappa
\nabla \alpha\big)
\ni\fDAMprime(\alpha)
&&\text{in }Q,
\end{align}\end{subequations}
when we confine ourselves to \eq{damage-avoiding-[0,1]-constraints}.
Let us note that we introduce the auxiliary variable $v$ standing for velocity
and write, rather for later purposes in Sect.\,\ref{sect-damage-discretisation}
the 1st-order system instead of the 2nd-order one.
The mathematical treatment relies on the linearity of \eq{8-GSM1**-u-special} 
in terms of $u$ but, on the other hand, \eq{8-GSM1**-z-special} is nonlinear 
in terms of $e=e(u)$.

Rather for 
simplicity, we consider the scenarios \eq{damage-avoiding-[0,1]-constraints},
which now means that $\bbC'(0)=0$ and possibly \NEW{(in the first option
in \eq{damage-avoiding-[0,1]-constraints})}
  also $\bbC'(1)=0$.
We consider a nested sequence of finite-dimensional subspaces 
of $H^1(\Omega;\R^d)$ and $H^1(\Omega)$ indexed by $k\in\N$ whose union
is dense in these Banach spaces, and then
an $H^1$-conformal Galerkin approximation,
denoting the approximate solution thus obtained by $(u_k,\alpha_k)$. 
For simplicity, we assume that 
$u_0,v_0\in V_1\subset V_k\subset H^1(\Omega;\R^d)$ as used for the Galerkin 
approximation; in fact, a natural qualification $v_0\in L^2(\Omega;\R^d)$ 
would in general need an approximation $v_{0,k}\in V_k$ such that
$v_{0,k}\to v_0$ strongly in $L^2(\Omega;\R^d)$. 

We allow for a {\it complete damage}\index{damage!complete} in the elastic 
response, although a resting Stokes-type viscosity due to $D_0$ is needed for 
the following \NEW{assertion relying on the linearity of
  $\pl_e\varphi(\cdot,\alpha)$, i.e.\ on that $\varphi(\cdot,\alpha)$ is
  quadratic}:   

\begin{proposition}[Existence in the linear model]
\label{prop-1}
Let the ansatz \eq{damage-gradient} 
be considered, let also 
$\varrho,\kappa\in L^\infty(\Omega)$ with ${\rm ess\,inf}\,\varrho>0$
and ${\rm ess\,inf}\,\kappa>0$, 
$f\in L^1(I;L^2(\Omega;\R^d))$, $g\in L^2(\Sigma;\R^d))$, 
$u_0\in H^1(\Omega;\R^d)$, $v_0\in L^2(\Omega;\R^d)$, 
$\alpha_0\in H^1(\Omega)$ with $0\le\alpha\COL{_0}\le 1$ a.e.\ on
$\Omega$ be supposed, 
$\zeta:\R\to\R^+$ be convex and lower semicontinuous with 
$\zeta(\cdot)\ge\epsilon|\cdot|^2$ for some $\epsilon>0$, and let
\eq{damage-avoiding-[0,1]-constraints} hold, and let also 
\eq{damage-gradient-a} be considered with
\begin{subequations}\label{ass-damage-visco}
\begin{align}\label{ass-damage-visco1}
&\text{$\bbC\In C^1([0,1];\R^{(d\times d)^2})$ be symmetric 
positive-semidefinite valued,} 
\\&\label{ass-damage-visco2}
\text{$\bbD(\cdot)=\bbD_0+\taur\bbC(\cdot)$ with $\taur\ge0$ and $\bbD_0$ 
symmetric positive-definite.} 
\end{align}\end{subequations}
Then the Galerkin approximation $(u_k,\alpha_k)$ exists and,
for selected subsequences, we have 
\begin{subequations}\label{conv-damage-linear}
\begin{align}\label{conv-damage-linear-u}
&u_k\to u\
&&\text{weakly* in }\ H^1(I;H^1(\Omega;\R^d))
\,\cap\,W^{1,\infty}(I;L^2(\Omega;\R^d))\ \text{ and}
\\&&&\label{conv-damage-linear-u+}\text{strongly in }
L^2(I;H^1(\Omega;\R^d))
\,,\ \text{ and}
\\&\alpha_k\to\alpha&&\text{weakly* in }\ 
H^1(I;L^2(\Omega))\,\cap\,L^\infty(I;H^1(\Omega))\,\cap\,L^\infty(Q),
\end{align}\end{subequations}
and every such a limit $(u,\alpha)$ is a weak solution in the sense of
Definition~\ref{def1} with $r_{\text{\sc c}}^{}=0$.
Moreover, even $e(\DT u_k)\to e(\DT u)$ strongly in $L^2(Q;\Rsym)$.
\end{proposition}

\medskip\noindent{\it Proof.}
The apriori estimates 
in the spaces occurring in (\ref{conv-damage-linear}a,c) 
can be obtained by standard
energetic test by $\DT u_k$ and $\DT\alpha_k$, which leads to 
\eq{KV-damage-engr} written for the Galerkin approximation,
and using H\"older's, Young's, and Gronwall's inequalities.

After selecting a  subsequence weakly* converging in the sense 
(\ref{conv-damage-linear}a,c) and using the 
Aubin-Lions theorem for the damage and then continuity of the
superposition operator induced by $\bbC(\cdot)$, we can 
pass to the limit first in the semilinear force-equilibrium
equation. We put $w:=u+\taur\DT u$ and write the
limit equation \eqref{8-GSM1**-u} as 
\begin{align}\label{KV-visco-damageble-strong--}
\frac{\varrho}{\taur}\DT w
-{\rm div}\big(\bbD_0e(\DT u)+
\bbC(\alpha)e(w)\big)=f+\frac{\varrho}{\taur}\DT u
\end{align}
accompanied with the corresponding initial/boundary conditions from 
\eq{8-GSM1**BC}--\eq{8-GSM1**IC}.

For the damage flow rule, we need the strong convergence
of $\{e(u_k)\}_{k\in\N}$, however. 
Furthermore, we denote $w_k:=u_k+\taur\DT u_k$ 
and, using the linearity of $\pl_e\varphi(\cdot,\alpha)$, write the 
Galerkin approximation of the force equilibrium as\footnote{\COL{More
    precisely, \eq{KV-visco-damageble-strong-} is to be understood
valued in $V_k^*$.}}
\begin{align}\label{KV-visco-damageble-strong-}
\frac{\varrho}{\taur}\DT w_k-{\rm div}\big(\bbD_0e(\DT u_k)+
\bbC(\alpha_k)e(w_k)\big)=f+\frac{\varrho}{\taur}\DT u_k\,.
\end{align}
Then we subtract \eq{KV-visco-damageble-strong--} and 
\eq{KV-visco-damageble-strong-}, and test it by
    $w_k-w$,
and 
integrate over the time interval $[0,t]$. This gives
\begin{align}\nonumber
&\int_\Omega\frac{\varrho}{2\taur}|w_k(T){-}w(T)|^2
+\frac12\bbD_0e(u_k(T){-}u(T)){:}e(u_k(T){-}u(T))\,\d x
\\[-.4em]&\nonumber\qquad
+\int_Q\bbD_0\taur e(\DT u_k{-}\DT u)\Colon e(\DT u_k{-}\DT u)
+\bbC(\alpha_k)e(w_k{-}w)\Colon e(w_k{-}w)\,\d x\d t
\\[-.4em]&
=\int_Q\big(\bbC(\alpha_k){-}\bbC(\alpha)\big)e(w)\Colon e(w_k{-}w)+
\frac{\varrho}{\taur}(\DT u_k{-}\DT u)
\Cdot(w_k{-}w)\,\d x\d t
\to0\,.
\label{KV-visco-damageble-conv-strong}
\end{align}
Here we used that $\DT u_k{-}\DT u\to0$ strongly in  $L^2(Q;\R^d)$ by the 
Aubin-Lions theorem and also that $(\bbC(\alpha_k){-}\bbC(\alpha))e(w)\to0$ 
strongly in $L^2(Q;\Rsym)$. This gives \eq{conv-damage-linear-u+}.
In fact, \eq{KV-visco-damageble-conv-strong} is again \NEW{a}
rather conceptual strategy
and still a strong approximation of $(u,w)$ is needed to facilitate 
usage of the Galerkin identity and convergence-to-zero of the additional  
terms thus arising.

The limit passage in the damage variational inequality towards 
\eq{id6-} is then easy by (semi)continuity.
\hfill$\Box$\medskip

In some applications a non-quadratic $\varphi(\cdot,\alpha)$ is a reasonable
ansatz in particular \COL{because} damage may act very differently on
compression than on tension, cf.\ \eqref{KV-damage-exa2} below.
Examples are concrete- or masonry-, or rock-type materials where mere 
compression practically does not cause damage while tension (as well as 
shear) may cause damage relatively easily.
Unfortunately, Proposition~\ref{prop-1} does not cover
such models. 
Two options allowing for $\alpha$-dependent $\bbD$ are doable: a 
unidirectional
damage \NEW{with hardening-like effect} and bi-directional (i.e.\ with possible
healing) damage.
Note that \eqref{ass-damage-visco} is not needed.
\NEW{In the first option, the constraint $\alpha\ge0$ is never active
 and $\alpha\le1$ is only ``semi-active'', both leading to zero 
Lagrange multiplier $r_{\text{\sc c}}^{}$.}

\begin{proposition}[Unidirectional damage 
in nonlinear models]\label{prop-2}
\index{damage!uniqueness}
Let 
the data $\varrho$, $\zeta(\cdot)$, $\kappa$, $f$, $g$, $u_0$, $v_0$, and 
$\alpha_0$ be as in Proposition~\ref{prop-1}, and let also 
$\varrho\in W^{1,\COL{r}}(\Omega)$ \COL{with $r=3$ for $d=3$ or
  $r>1$  for $d=2$},
\begin{align}\label{KV-damage-ass-phi+weaker}
|\pl_e\varphi(e,\alpha)|\le C\big(1+|e|\big)\ \ \ \text{ and }\ \ \ 
|\pl_\alpha\varphi(e,\alpha)|\le C\big(1+|e|^2\big).
\end{align}
Let moreover 
 $\bbD:[0,1]\to\R^{(d\times d)^2}$ be symmetric-valued, continuous, 
monotone (nondecreasing) with respect
to the L\"owner ordering (i.e.\ ordering of $\Rsym$ by the cone of positive 
semidefinite matrices), and with $\bbD(0)$ positive definite.
Let moreover the second option in \eq{damage-avoiding-[0,1]-constraints} hold.
%
Then the 
Galerkin approximate solutions do exist with $r_{\text{\sc c},k}^{}=0$.
The sequence $\{(u_k,\alpha_k)\}_{k\in\N}$ possesses subsequences 
such that again \eqref{conv-damage-linear}
hold. 
The limit of each such a subsequence solves the initial-boundary-value problem 
\eq{8-GSM1**}--\eq{8-GSM1**IC}
weakly in the sense of
Definition~\ref{def1} \NEW{with $r_{\text{\sc c}}^{}=0$}. 
%
\end{proposition}

\medskip\noindent{\it Proof.}
We perform the test of (\ref{8-GSM1**}a,b) in the Galerkin 
approximation by $(\DT u_k,\DT\alpha_k)$. By using 
the data qualification and H\"older and Gronwall inequalities, this gives the 
estimates in the spaces occurring in (\ref{conv-damage-linear}a,c). 
By comparison from \eq{8-GSM1**-u}, we obtain also the bound\footnote{More 
precisely, this bound is valid only in Galerkin-induced seminorms or 
for the Hahn-Banach extension, cf.\ \cite[Sect.\,8.4]{Roub13NPDE}.} for 
$\DDT u_k$ \COMMENT{CHECK Book4 FOR $\varrho\DDT u_k$!!} in $L^2(I;H^1(\Omega;\R^d)^*)$ by estimating
\begin{align}\nonumber
&\int_Q\DDT u_kv\,\d x\d t
=\int_Q(f+{\rm div}\,\sigma_k)\Cdot\frac v\varrho\,\d x\d t
=\int_Q\frac{f\Cdot v}\varrho
-\sigma_k\Colon\nabla\frac v\varrho
\,\d x\d t
\\&\ \ =\int_Q\frac{f\Cdot v}\varrho
-\frac{\sigma_k\Colon\nabla v}\varrho
+\frac{\sigma_k\Colon\nabla\varrho}{\varrho^2}\,\d x\d t
+\int_\Sigma\frac{g\Cdot v}\varrho\,\d S\d t
\le C\|v\|_{L^2(I;H^1(\Omega;\R^d))},
\label{est-of-DDTu}
\end{align}
where $C$ is dependent on the already obtained estimates 
(\ref{conv-damage-linear}a,c); note also that we need a certain 
smoothness of $\varrho$, as supposed.


After selection of weakly* convergent subsequences, 
we prove the strong convergence 
%
\eqref{conv-damage-linear-u+}.
We use a slightly different estimation comparing to 
\eqref{KV-visco-damageble-conv-strong} based on 
a test by $\DT u_k{-}\DT u$, namely now
we employ the test by $u_k{-}u$ to estimate
\begin{align}\nonumber
&\int_\Omega\frac12\bbD(\alpha_k(t)) e(u_k(t){-}u(t))\Colon e(u_k(t){-}u(t))
\,\d x
+\int_0^t\!\!\int_\Omega\big(\pl_e^{}\varphi(e(u_k),\alpha_k)-\pl_e^{}\varphi(e(u),\alpha_k)\big)
\Colon e(u_k{-}u)\,\d x\d t
\\&\nonumber=\int_0^t\!\!\int_\Omega\frac{\pl}{\pl t}
\frac12\bbD(\alpha_k)e(u_k{-}u)\Colon e(u_k{-}u)
+\big(\pl_e^{}\varphi(e(u_k),\alpha_k)-\pl_e^{}\varphi(e(u),\alpha_k)\big)
\Colon e(u_k{-}u)\,\d x\d t
\\&\nonumber=
\int_0^t\!\!\int_\Omega
\big(\bbD(\alpha_k) e(\DT u_k{-}\DT u)+\pl_e^{}\varphi(e(u_k),\alpha_k)
-\pl_e^{}\varphi(e(u),\alpha_k)\big)
\Colon e(u_k{-}u)
+\frac12\DT\alpha_k\bbD'(\alpha_k)e(u_k{-}u)\Colon e(u_k{-}u)\,\d x\d t
\\&\nonumber\le\int_0^t\!\!\int_\Omega(f-\varrho\DDT u_k)\Cdot(u_k{-}u)
-\big(\bbD(\alpha_k)e(\DT u)+\pl_e^{}\varphi(e(u),\alpha_k)\big)\Colon e(u_k{-}u)
\,\d x\d t
\\&\nonumber=
\int_0^t\!\!\int_\Omega f\Cdot(u_k{-}u)+\varrho\DT u_k\Cdot(\DT u_k{-}\DT u)
-\big(\bbD(\alpha_k) e(\DT u)
+\pl_e^{}\varphi(e(u),\alpha)\big)\Colon e(u_k{-}u)\,\d x\d t
\\[-.4em]&
\qquad\qquad
+\int_0^t\!\!\int_\Omega
\big(\pl_e^{}\varphi(e(u),\alpha){-}\pl_e^{}\varphi(e(u),\alpha_k)\big)
\Colon e(u_k{-}u)\,\d x\d t
%
-\int_\Omega\varrho\DT u_k(t)\Cdot(u_k(t){-}u(t))\,\d x\ \to\ 0.
\label{KV-damage-strong-e(t)}
\end{align}
because $\DT\alpha_k\bbD'(\alpha_k)e(u_k{-}u)\Colon e(u_k{-}u)\le0$
a.e.\ on $Q$ since $\DT\alpha_k\le0$ due to the assumption that
$\zetai(\DT\alpha)=+\infty$ for $\DT\alpha>0$ and $\bbD'(\cdot)$ is positive
semidefinite due to the assumption of monotone dependence of $\bbD(\cdot)$.

Then, having this strong convergence, we can easily pass to the limit by 
(semi)con\-tinuity both towards the identity \eq{7-very-weak-sln2} and towards
variational inequality \eq{id6-}.
\hfill$\Box$\medskip

It is interesting that the usual ``limsup-argument'' relying on the
energy conservation to prove the strong convergence \eq{conv-damage-linear-u+} 
could not be used while \eq{KV-damage-strong-e(t)} worked, relying 
on the unidirectionality of damage evolution. Let us 
further illustrate the opposite situation when \eq{KV-damage-strong-e(t)} 
does not work while the energy conservation holds
and facilitates the mentioned limsup-argument:

\begin{proposition}[Damage with healing in nonlinear models]\label{prop-2+}
\index{damage!uniqueness}
Let the data $\varrho$, $\bbD(\cdot)$, $\kappa$, $f$, $g$, 
$u_0$, $v_0$, and $\alpha_0$ be as in Proposition~\ref{prop-2}, and again
\eq{KV-damage-ass-phi+weaker} hold\COL{.}
Let moreover 
\begin{subequations}\label{option-b}\begin{align}
&\nonumber
\exists\,0<\eps\le C\ \forall (e,z)\in\Rsym\times\R:\ \ \ 
\eps|z|^2\le\zeta(z)\le C(1{+}|z|^2),\ \ \
\\&\hspace{14.2em}
|\pl_\alpha\varphi(e,z)|\,\le\, C(1+|e|).
\label{option-b2}\end{align}\end{subequations}
and the Galerkin approximation is $H^2$-conformal
so that ${\rm div}(\kappa\nabla\alpha_k)$ 
\NEW{is well defined}.
Then the mentioned Galerkin approximate solutions do exist.
The sequence $\{(u_k,\alpha_k,r_{\text{\sc c},\COL{k}}^{})\}_{k\in\N}$ possesses subsequences 
such that again \eqref{conv-damage-linear}
hold together with 
\begin{align}\label{r-conv}
r_{\text{\sc c},\COL{k}}^{}\to r_{\text{\sc c}}^{}\ \ \text{ weakly in }\ L^2(Q).
\end{align}
The limit of each such a subsequence solves the initial-boundary-value problem 
\eq{8-GSM1**}--\eq{8-GSM1**IC}
weakly in the sense of
Definition~\ref{def1}. Moreover,
${\rm div}(\kappa\nabla\alpha)\in L^2(Q)$, the damage flow rule
(\ref{8-GSM1**}b,c) holds a.e.\ on $Q$, and the energy conservation holds.

\end{proposition}

\medskip\noindent{\it Proof.}
Let us outline only the differences from the proof of Proposition\ref{prop-2+}.
Beside the a-priori estimates there, we further test the approximated damage 
flow-rule by ${\rm div}(\kappa\nabla\alpha_k)$
We thus obtain a bound for ${\rm div}(\kappa\nabla\alpha_k)$ in 
$L^2(Q)$, and eventually also for 
$r_{\text{\sc c},k}^{}\in
{\rm div}(\kappa\nabla\alpha_k)+\fDAMprime(\alpha_k)
-\pl\zeta(\DT\alpha_k)-\pl_\alpha\varphi(e(u_k),\alpha_k)
$ 
in $L^2(Q)$.\footnote{Here we have employed also the calculus
$\int_\Omega r_{\text{\sc c},k}^{}
{\rm div}\big(\kappa\nabla\alpha_k\big)\,\d x=
\int_\Omega\pl\delta_{[0,1]}^{}(\alpha_k)
{\rm div}\big(\kappa\nabla\alpha_k\big)\,\d x
=-\int_\Omega\kappa\nabla\big[\pl\delta_{[0,1]}^{}(\alpha_k)\big]
\Cdot\nabla\alpha_k\,\d x
=-\int_\Omega\kappa\pl^2\delta_{[0,1]}^{}(\alpha_k)
\nabla\alpha_k\Cdot\nabla\alpha_k\,\d x\le0$.
}
 Here we used also the growth conditions \eqref{option-b2}
guaranteeing that both $\pl\zeta(\DT\alpha_k)$ and 
$\pl_\alpha\varphi(e(u_k),\alpha_k)$ are bounded in $L^2(Q)$. 

We now can pass to the limit in the force equilibrium just \COL{by} the 
weak convergence and monotonicity of $\pl_e\varphi(\cdot,\alpha_k)$
and the Aubin-Lions compactness theorem used for $\alpha_k$. 
Having the limit equation \eq{8-GSM1**-u} at disposal in the 
weak sense \eq{7-very-weak-sln2}, we can test it by $v=\DT u$ and 
show energy conservation in this part of the system. To this goal,
it is important that both $\DT\alpha$ and 
$\pl_\alpha\varphi(e(u),\alpha)$ are in $L^2(Q)$ so that the 
chain rule \eq{chain-rule-for-GSM} rigorously holds 
and that $\sqrt\varrho\DDT u\in L^2(I;H^1(\Omega;\R^d)^*)$ 
is in duality with $\sqrt\varrho\DT u\in L^2(I;H^1(\Omega;\R^d))$
so that also the chain rule $\int_Q\varrho\DT u\Colon\DDT u\,\d x\d t
=\int_\Omega\frac12\varrho|\DT u(T)|^2-\frac12\varrho|\DT u(0)|^2\,\d x$;
the information about $\sqrt\varrho\DDT u$ can be obtained by a
simple modification of \eq{est-of-DDTu}. By this test, we obtain
\begin{align}\nonumber
&\hspace{0em}
\int_\Omega\frac12\varrho|\DT u(T)|^2+\varphi(e(u(T)),\alpha(T))\,\d x
+\int_Q\!\bbD(\alpha)e(\DT u)\Colon e(\DT u)\,\d x\d t=\!\int_{\Gamma}\!g\Cdot\DT u\,\d S\d t
\\[-.4em]&
\qquad\qquad
+\int_\Omega\frac12\varrho|v_0|^2+\varphi(e(u_0),\alpha_0)\,\d x
+\!\int_Qf\Cdot\DT u-\DT\alpha\pl_\alpha\varphi(e(u),\alpha)\,\d x\d t
\label{7-very-weak-sln2-tested}
\end{align}
Instead of \eq{KV-damage-strong-e(t)}, we now estimate by weak semicontinuity
\begin{align}\nonumber
&\int_Q\bbD(\alpha)e(\DT u)\Colon e(\DT u)\,\d x\d t
  \le\liminf_{k\to\infty}\NEW{\int_Q}
  \bbD(\alpha_k)e(\DT u_k)\Colon e(\DT u_k)\,\d x\d t
\\&\nonumber
\quad=\int_\Omega\frac12\varrho|v_0|^2+\varphi(e(u_0),\alpha_0)\,\d x
+\lim_{k\to\infty}\!\int_Qf\Cdot\DT u_k-\DT\alpha_k\pl_\alpha\varphi(e(u_k),\alpha_k)\,\d x\d t
\\[-.4em]&\nonumber
\qquad\qquad-\liminf_{k\to\infty}
\int_\Omega\frac12\varrho|\DT u_k(T)|^2+\varphi(e(u_k(T)),\alpha_k(T))\,\d x
\\&\nonumber
\quad\le\int_\Omega\frac12\varrho|v_0|^2+\varphi(e(u_0),\alpha_0)\,\d x
-\!\int_\Omega\frac12\varrho|\DT u(T)|^2+\varphi(e(u(T)),\alpha(T))\,\d x
\\[-.3em]&\qquad\qquad+\!\int_Qf\Cdot\DT u-\DT\alpha\pl_\alpha\varphi(e(u),\alpha)\,\d x\d t
=\int_Q\bbD(\alpha)e(\DT u)\Colon e(\DT u)\,\d x\d t
\end{align}
where the last equality is due to \eq{7-very-weak-sln2-tested}. Altogether,
we \COL{have} proved that
$\liminf_{k\to\infty}\NEW{\int_Q}\bbD(\alpha_k)e(\DT u_k)\Colon e(\DT u_k)\,\d x\d t
=\int_Q\bbD(\alpha)e(\DT u)\Colon e(\DT u)\,\d x\d t$. From this,
we obtain even $e(\DT u_k)\to e(\DT u)$ strongly in $L^2(Q;\Rsym)$. 
More in detail, using uniform positive definiteness of $\bbD(\cdot)$, we 
perform the estimate
\begin{align}\nonumber&\!\!\!\!\!
\min_{\alpha\in[0,1]}|\bbD^{-1}(\alpha)|^{-1}\|e(\DT u_k{-}\DT u)\|_{L^2(Q;\R^{d\times d})}^2
\le\int_Q\bbD(\alpha_k)e(\DT u_k{-}\DT u)\Colon e(\DT u_k{-}\DT u)\,\d x\d t
\\&\!\!=\int_Q\bbD(\alpha_k)e(\DT u_k)\Colon e(\DT u_k)
-2\bbD(\alpha_k)e(\DT u_k)\Colon e(\DT u)+\bbD(\alpha_k)e(\DT u)\Colon e(\DT u)\,\d x\d t\to0.\!\!
\end{align}
Hence, we obtained even more that the desired strong convergence 
$e(u_k)\to e(u)$.

Now, beside the limit passage as in the proof of
\NEW{Proposition~\ref{prop-2+}},
also the limit passage towards the inclusion \eq{8-GSM1**-r} using that 
$N_{[0,1]}(\cdot)$ has a closed monotone graph is easy since 
$\alpha_k\to\alpha$ strongly in $L^2(Q)$ due to Aubin-Lions theorem while
$r_{\text{\sc c},k}^{}\to r_{\text{\sc c}}^{}$ weakly in $L^2(Q)$.
In particular, we have the chain rule
$\int_Qr_{\text{\sc c}}^{}\DT\alpha\,\d x\d t=0$ at disposal.

Eventually, as both $\DT\alpha$ and ${\rm div}(\kappa\nabla\alpha)$ 
are in $L^2(Q)$, we have also the chain rule
$\int_Q\DT\alpha{\rm div}(\kappa\nabla\alpha)\,\d x\d t=
\int_\Omega\frac12\kappa|\nabla\alpha_0|^2-\frac12\kappa|\nabla\alpha(T)|^2\,\d x$
and we can test the damage flow rule by $\DT\alpha$, and then sum it with 
\eq{7-very-weak-sln2-tested} to obtain the energy balance \eq{KV-damage-engr}.
\hfill$\Box$\medskip


Let us notice that a unidirectional damage \NEW{eveolution}
(i.e.\ $\zeta(\DT\alpha)=+\infty$ 
for $\DT\alpha>0$) together with the
damage, where the indicator 
function $\delta_{[0,1]}$ in \eqref{damage-gradient-E} and \eqref{8-GSM1**-r} 
must be considered, is not covered by Propositions~\ref{prop-1}--\ref{prop-2+}.
A particular case when $\zeta(\cdot)$ is positively homogeneous 
(i.e.\ the damage-process itself is rate-independent) allows a particular
treatment by using a so-called energetic formulation, invented for 
rate-independent systems by A.\,Mielke at al.\ 
\cite{Miel05ERIS,MieRou15RIST,MieThe04RIHM},
and later adapted for dynamical systems containing rate-independent sub-systems 
in \cite{Roub09RIPV}. We need a space of functions $I\to L^1(\Omega)$ 
of bounded variations, denoted by ${\rm BV}(I;L^1(\Omega))$,
i.e.\ the Banach space of functions $\alpha:I\to L^1(\Omega)$ 
with $\sup_{0\le t_0<t_1<...t_N\le T,\ N\in\N}\sum_{i=1}^N\|\alpha(t_i)-\alpha(t_{i-1})\|_{L^1(\Omega)}$ finite.

\begin{definition}[Energetic formulation]\label{def2}
A pair $(u,\alpha)\in H^1(I;H^1(\Omega;\R^d))\times {\rm BV}(I;L^1(\Omega))
\cap L^\infty(I;H^1(\Omega))$ is called an energetic solution to the 
initial-boundary-value problem \eq{8-GSM1**}--\eq{8-GSM1**IC}
if again \eq{7-very-weak-sln2} holds 
for all $v\In L^2(I;H^1(\Omega;\R^d))\,\cap\,H^{2}(I;L^2(\Omega;\R^d))$
with $v|_{t=T}^{}=\DT{v}|_{t=T}^{}=0$, if also the energy balance 
\eq{KV-damage-engr} holds, and if the so-called semi-stability
\begin{align}&\nonumber
\int_\Omega\varphi\big(e(u(t)),\alpha(t)\big)+\frac\kappa2|\nabla\alpha(t)|^2
\le\int_\Omega\varphi\big(e(u(t)),z\big)+\frac\kappa2|\nabla z|^2
+\zeta(z{-}\alpha(t)\big)\d x
\end{align}
holds for all $t\in I$ and all $z\in H^1(\Omega)$ with $0\le z\le1$ 
a.e.\ on $\Omega$, and also $\alpha|_{t=0}=\alpha_0$.
\end{definition}

Under the semistability of the initial damage profile $\alpha_0$, existence of 
energetic solutions can be shown by the implicit or semi-implicit 
time discretisation and by an explicit construction of a so-called 
mutual recovery sequence \cite{MiRoSt08GLRR}.\footnote{This specifically
here means that, having $u_k\to u$ strongly in $H^1(\Omega;\R^d)$ and 
$\alpha_k\to\alpha$ weakly in $H^1(\Omega)$ with $0\le\alpha_k\le1$, we need 
to find a sequence $\{z_k\}_{k\in\N}$ such that 
\begin{align*}
\limsup_{k\to\infty}
\int_\Omega\varphi\big(e(u_k),z_k\big)+\frac\kappa2|\nabla z_k|^2
+\zeta(z_k{-}\alpha_k\big)-\varphi\big(e(u_k),\alpha_k\big)
-\frac\kappa2|\nabla\alpha_k|^2\d x\qquad
\\[-.4em]\le\int_\Omega\varphi\big(e(u),z\big)+\frac\kappa2|\nabla z|^2
+\zeta(z{-}\alpha\big)-\varphi\big(e(u),\alpha(t)\big)
-\frac\kappa2|\NEW{\nabla\alpha}|^2\d x\,.
\end{align*} 
This needs quite sophisticated construction devised in 
\cite{ThoMie10DNEM}. When the quadratic gradient term in 
\eq{damage-gradient-E} would be replaced by the $p$-power 
with $p>d$, a simpler construction would apply, cf.\ \cite{ThoMie10DNEM}.}
\COL{For even an anisothermally enhanced model, we also refer to
\cite{LRTT18RIDT}.}

An important special case consist in isotropic materials, where one
can easily distinguish response under volumetric and shear load which
might be very diverse. To this goal, we  use the 
decomposition of the strain to its compression/tension spherical and 
its deviatoric parts:
\begin{align}\label{decomposition-of-strain}
e={\rm sph}^+e+{\rm sph}^-e+{\rm dev}\,e\ \ \ \text{ with }\ \ \ 
{\rm sph}^\pm\,e=\frac{({\rm tr}\,e)^\pm}d\bbI\,,
\end{align}
where $({\rm tr}\,e)^+=\max(\sum_{i=1}^de_{ii},0)$ and 
$({\rm tr}\,e)^-=\min(\sum_{i=1}^de_{ii},0)$.
Note that the deviatoric and the spherical strains 
from \eq{decomposition-of-strain} are orthogonal to each other
and the deviatoric strain is trace-free, i.e.\ ${\rm dev}\,e\Colon
{\rm sph}^\pm\,e=0$ and ${\rm tr}({\rm dev}\,e)=0$. This decomposition 
allows for distinguishing the response under compression (usually not causing 
damage), tension (so-called Mode I damage), or shearing (Mode II damage).
A combination is called a mixed mode, and 
altogether we speak about a mode sensitive damage.

The latter growth restriction in \eq{option-b} excludes
quadratic energies of the type \eqref{damage-gradient-a}.
An example of a model with the mode-sensitive isotropic stored energy 
satisfying this restriction 
and the damage dissipation potential with quadratic coercivity might be
\begin{subequations}\label{mixity-model}\begin{align}
&\label{KV-damage-exa2}
\varphi(e,\alpha)=
\frac d2\bigg(K(1)|{\rm sph}^-\!e|^2+
\frac{K(\alpha)|{\rm sph}^+\!e|^2}{\sqrt{1{+}\epsilon({\rm tr}\,e)^2}}
\bigg)
+\frac{G(\alpha)|{\rm dev}\,e|^2}{\sqrt{1{+}\epsilon|{\rm dev}\,e|^2}}\,,
\\&\label{zeta}
\zeta(e;\DT\alpha)=\begin{cases}-\Frakg_{\rm c}\DT\alpha
+\nu\DT\alpha^2&\text{if $\DT\alpha\le0$},\\
+\infty&\text{if $\DT\alpha>0$},\end{cases}
\end{align}\end{subequations}
where $K=K(\alpha)$ is the bulk modulus and $G=G(\alpha)$ is the shear modulus;
recall that $K=\lambda+2G/d$ with $\lambda$ and $G$ the so-called Lam\'e 
constants. The coefficient $\Frakg_{\rm c}>0$ in \eq{zeta} is called a fracture
toughness while the coefficient $\nu>0$ makes fast damage more dissipative
(more heat producing) than slower damage, which might be sometimes relevant
and which makes mathematics sometimes easier, as in Propositions~\ref{prop-1} 
and \ref{prop-2} above. 
The (small) regularizing parameter $\epsilon>0$ makes 
the tension and shear stress bounded if $|e|$ is (very) large and makes
the growth restriction on $\pl_\alpha\varphi$ in \eq{option-b}
satisfied, while $\eps=0$ is admitted \NEW{in the case of}
the second option in \eq{damage-avoiding-[0,1]-constraints}.

Let us illustrate heuristically how the flow-rule 
$\zeta(\DT\alpha)+\pl_{\alpha}\varphi(e,\alpha)\ni0$ with the
initial condition $0<\alpha(0)=\alpha_0\le1$ operates when the 
loading gradually increases.
For the example \eq{KV-damage-exa2} with $\epsilon=0$ and \eq{zeta} with 
$\nu=0$, the stress is $\sigma=\pl_e\varphi(e,\alpha)=
\sigma_{\rm sph}^-+\sigma_{\rm sph}^++\sigma_{\rm dev}^{}$
with $\sigma_{\rm sph}^-=\NEW{dK(1)
{\rm sph}^-e}$, $\sigma_{\rm sph}^+=dK(\alpha){\rm sph}^+e$, and $\sigma_{\rm dev}=
2\SHEAR(\alpha){\rm dev}\,e$.
The driving force for damage evolution expressed in therm of 
the actual stress is
$$
\pl_{\alpha}\varphi(e,\alpha)=\frac12K'(\alpha)|{\rm sph}^+e|^2+
G'(\alpha)|{\rm dev}\,e|^2
=\frac{K'(\alpha)}{2dK(\alpha)^2}|\sigma_{\rm sph}^+|^2+
\frac{G'(\alpha)}{4\SHEAR(\alpha)^2}|\sigma_{\rm dev}|^2.
$$
Then the criterion $\pl_{\alpha}\varphi(e,\alpha_0)=\Frakg_{\rm c}$
reveals the stress needed to start damaging the material. In the pure 
shear or pure tension, this critical stress is
\index{stress!effective fracture} 
\begin{align}\label{KV-damage-exa1-toughness}
&|\sigma_{\rm dev}|=\SHEAR(\alpha_0)\sqrt{\frac{4\Frakg_{\rm c}}{\SHEAR'(\alpha_0)}}
=\text{``effective fracture stress'' in Mode II}.
\\&|\sigma_{\rm sph}^+|=K(\alpha_0)\sqrt{\frac{2d\Frakg_{\rm c}}{K'(\alpha_0)}}
=\text{``effective fracture stress'' in Mode I}.
\end{align}
respectively.
If $G(\cdot)G'(\cdot)^{-1/2}$ and $K(\cdot)K'(\cdot)^{-1/2}$ are increasing 
(in particular if $G(\cdot)$ and $K(\cdot)$ are concave), and the loading is 
via stress rather than displacement, damage then accelerates when started so 
that the rupture happens immediately (if any rate and spatial-gradient effects 
are neglected). 


\section{Phase-field concept towards fracture}\label{sect-damage-PF-fracture}

The concept of bulk damage can (asymptotically) imitate the philosophy of 
fracture along surfaces (\emph{cracks})\index{crack} provided the damage 
stored energy $\fDAM$ is big. A popular ansatz takes the basic model 
\eq{damage-gradient} with \eq{damage-gradient-a} for
\begin{align}
\bbC(\alpha):=\Big(
\frac{\eps^2}{\eps_0^2}
{+}\alpha^2\Big)\bbC_1,\ \ \ \ 
\fDAM(\alpha):=-\Frakg_{\rm c}
\frac{(1{-}\alpha)^2}{2\eps},\quad\text{ and }\quad 
\kappa:=\eps \Frakg_{\rm c}
\end{align}
with $\Frakg_{\rm c}$ denoting the energy of fracture
and with $\eps$ controlling a ``characteristic'' width of the 
{\it phase-field fracture}\index{fracture!phase-field approach} zone(s). 
This width is supposed to be small with respect to the size of the whole body. 
Then, \eq{damage-gradient-E} looks (up to the forcing $f$ and $g$) as
\begin{align}
&\!\!\!\!\!\!
\mathscr{E}(u,\alpha):=\int_\Omega
\COL{\gamma(\alpha)}
\bbC e(u){:}e(u)
+
\!\!\lineunder{
\Frakg_{\rm c}\Big(\frac{1}{2\eps}(1{-}\alpha)^2\!
+\frac\eps2|\nabla\alpha|^2\Big)}{crack surface density}\!\,\d x
\COL{\ \ \text{ with }\ \!\!\lineunder{\gamma(\alpha)=\frac{
(\eps/\eps_0)^2{+}\alpha^2}2}{degradation function}}
\label{eq6:AM-engr+}
\end{align}
\COL{and } with $\eps_0>0$. The physical dimension of $\eps_0$ as well as of 
$\eps$ is m (meters) while the physical dimension of $\Frakg_{\rm c}$ 
is J/m$^2$. This is known as the so-called \emph{Ambrosio-Tortorelli 
functional}\index{Ambrosio-Tortorelli functional}.
\footnote{In the static case, 
this approximation was proposed by Ambrosio and Tortorelli 
\cite{AmbTor90AFDJ,AmbTor90AFDP} \NEW{originally}
for the scalar Mumford-Shah functional \cite{MumSha89OAPS} and the asymptotic 
analysis for $\eps\to0$ was rigorously executed.
A
generalization (in some sense in the spirit of finite-fracture 
mechanics) is in \cite{conti.focardi.iurlano2016}.
The generalization for the vectorial case is in \cite{Foca01VAFD,focardi.iurlano2014,iurlano2014}.
Later, it was extended for evolution situation, namely for 
a rate-independent
damage, in \cite{Giac05ATAQ}, see also 
also \cite{BoFrMa08VAF,BoLaRi11TDMD,caponi2018,LaOrSu10ESRM,MieRou15RIST}
where also inertial forces are sometimes considered.}
The fracture toughness $\Frakg_{\rm c}$ is now involved in 
\eq{eq6:AM-engr+} \NEW{instead of}
  the dissipation potential \eq{zeta}\NEW{, i.e.\ $\zeta$ in \eq{zeta}}
is now considered with $\Frakg_{\rm c}=0$.

It should be emphasized that, in the ``crack limit'' for $\eps\to0$, 
the phase-field fracture model 
\eq{eq6:AM-engr+} approximates (at least in the static and quasistatic cases) 
the 
true infinitesimally thin cracks
in Griffith's \cite{Grif21PRFS} variant (i.e.\ competition 
of energies), which works realistically for crack propagation 
but might have unrealistic difficulties with crack initiation,\footnote{In fact,
as $\fDAMprime(1)=0$, the initiation of damage has zero threshold and is 
happening even on very low stress but then, if $\eps>0$ is very small, stops 
and high stress is needed to continue damaging.} while scaling 
of the fracture energy to 0 if $\eps\to0$ might lead to opposite effects,
cf.\ also the discussion e.g.\ in Remark~\ref{rem-stress-vs-engr} below.
This is partly reflected by \NEW{the fact} that, in its rate-independent
variant, the damage 
and phase-field fracture models admit many various solutions of very 
different characters, as presented in \cite{MieRou15RIST}.
In the dynamic variant, the influence of overall stored energy during
fast rupture (which may be taken into account during quasistatic evolution) 
seems eliminated because of finite speed of propagation of information
about it.\footnote{Yet, in this dynamic case, the analysis for $\eps\to0$
remains open and, even worse, in the limit crack problem one should 
ca\COL{r}e about non-interpenetration, which is likely very difficult; cf.\ the
analysis for the damage-to-delamination problem \cite{MiRoTh12DDNE}.}

Various modifications have been devised. For example,  
Bourdin at al.\ \cite{BMMS14MPCC} used 
\begin{align}\label{BMMS}
\fDAM(\alpha)=\frac{3\Frakg_{\rm c}\alpha}{8\eps}\ \ \text{ 
and }\ \ \kappa=\frac{3\Frakg_{\rm c}\eps}4
\end{align} 
in \eq{damage-gradient} 
and $\Frakg_{\rm c}$ the energy of fracture in J/m$^2$ 
and with $\eps$ controlling a ``characteristic'' width of the 
phase-field fracture zone. 
Although it activates damage process only when stress achieves some threshold,
it 
exhibits a similar undesired behaviour as \eq{eq6:AM-engr+} when $\eps\to0$ 
and leads to consider $\eps>0$ as another parameter (without intention
to put it 0) in addition to $\Frakg_{\rm c}$ to tune the model.
\COL{Various modifications of the degradation function from
\eq{eq6:AM-engr+} therefore appeared in literature. E.g.\
a cubic degradation function $\gamma$ has been used e.g.\ in
\cite{Bord12IAPF,VMBV14PFMB}.}
Inspired by \eq{KV-damage-exa1-toughness}, keeping still the 
original motivation $\eps\to0$, one can think about 
some
$\gamma$ convex increasing with $\gamma'(1)=\mathscr{O}(1/\eps)$.
\COL{The mentioned cubic ansatz is not compatible with
  these requirement. Some more sophisticated $\gamma$'s have been
  devised in \cite{SKBN18HAPF,WuNgu18LSIP}.}
%
We have thus
\COL{more}
independent parameters \COL{than only $\Frakg_{\rm c}$ and $\eps$ in
\eq{eq6:AM-engr+}} 
to  specify the \COL{lengh-scale of damage zone},
fracture propagation and \COL{fracture} initiation.

\begin{remark}[Finite fracture mechanics (FFM)]\label{rem-FFM}\upshape
\COL{In contrast to the Griffith model relevant rather for 
infinitesimally short increments of cracks}\footnote{\NEW{Cf.\ e.g.\ the
analysis and discussion in \cite{SicMar13GDLG} in the quasistatic
situations.}}\COL{, the finite (large) increments needs rather the concept of
energetic solution.}
As already mentioned, it does not seem much realistic to count with the 
overall strain energy in very distant spots (particularly in dynamical 
problems with finite speed of propagation of information), so rather 
only energy around a current point $x\In\Omega$ is to be considered
and cracks can propagate only by finite distance during incremental 
stepping. This is concept is commonly called a finite fracture mechanics 
(FFM); this term has been suggested by Z.\,Hashin \cite{Hash44FTFC}, but 
being developed rather gradually by several authors, see e.g.\ 
\cite{TaCoPu05FMFC}. It occurs useful in particular in quasistatic
problems \COL{which neglect inertia} to compensate 
(rather phenomenologically) \COL{this simplification}.
\end{remark}

\begin{remark}[Coupled stress-energy criterion]\label{rem-stress-vs-engr}\upshape
In addition to FFM, in fracture (or in general damage) 
mechanics, there is a disputation whether only sufficiently big 
stress can lead to rupture or (in reminiscence to Grifith's concept) whether 
(also or only) some sufficiently big energy in the specimen or around the 
crack process zone is needed for it. A certain standpoint is that both 
criteria should be taken into account.
This concept is nowadays referred as coupled stress-energy 
criterion.\footnote{Cf.\ 
the survey \cite{WeLeBe16RFFM},
and has been devised and implemented in many variants in engineering literature,
cf.\ e.g.\ \cite{CPCT06FFMC,GCIM16NSTC,Legu02STCC,Mant09ICOC,Mant14PIGC},
always without any analysis of numerical stability and convergence and thus
computational simulations based on these models, whatever practical 
applications 
they have, stay in the position of rather speculative playing with computers.}
Here, this coupled-criterion concept can be reflected by making $\zeta$ 
dependent on the strain energy. Having in mind FFM,\COMMENT{????} one can think to let 
$\zeta=\zeta(\widetilde\eps;\DT\alpha)$ with 
$\widetilde \eps(x)=\int_\Omega k(x{-}\widetilde x)
\varphi(e(u(\widetilde x)),\alpha(\widetilde x))\,\d\widetilde x$ 
for some kernel $k:\R^d\to\R^+$, making $-\pl_{\DT\alpha}\zeta(\widetilde
\eps;\DT\alpha)$ larger \COL{if} $\widetilde\eps$ is small. \COMMENT{....averadge strain 
energy density criterion  --- MAYBE:\\ 
Strain-energy-density factor applied to mixed mode crack problems
G.C. Sih - International Journal of fracture, 1974\\
OR\\
Application of an average strain energy density criterion to obtain the mixed 
mode fracture load of granite rock tested with the cracked asymmetric 
four-point bend .....
S.M.J. Razavi, M.R.M. Aliha, F. Berto - Theoretical and Applied Fracture..., 2017\\
OR\\
Control volumes and strain energy density under small and large scale yielding due to tension and torsion loading
P/ Lazzarin, F/ Berto - … \& Fracture of Engineering Materials ...., 2008 
}
\end{remark}

\begin{remark}[Mixity-mode sensitive cracks]\label{rem-mode-sensitive}\upshape
Combining the mixity-mode sensitive model \eq{KV-damage-exa2} with the crack 
surface density from \eq{eq6:AM-engr+}, one can distinguish the fracture by 
tension while mere compression does not lead to fracture, cf.\ 
\cite{MiWeHo10TCPF,SWKM14PFAD},
and one can also distinguish the {\it Mode I} 
({\it fracture}\index{fracture!mode I/II} by opening) from
{\it Mode II} (fracture by shear)\NEW{, cf.\ \cite{LanRoy09VAFM}}. 
\COL{More in detail, like in \eq{KV-damage-exa2}, we can use different
  degradation functions $\gamma$'s for the deviatoric part and the spherical
  compressive part and the spherical tension part.}
\end{remark}

\begin{remark}[Various other models]\upshape
An alternative option how to distinguish Mode I from Mode II is in the 
dissipation potential, reflecting the experimental observation that Mode II 
needs (dissipates) more energy than Mode I. Thus, one can take a 
state-dependent $\zeta=\zeta(e;\DT\alpha)$ e.g.\ as \eq{zeta} with 
$\Frakg_{\rm c}=\Frakg_{\rm c}({\rm sph}\,e,{\rm dev}\,e)>0$ to be rather small 
if ${\rm tr}e\gg|{\rm dev}\,e|$ (which indicates Mode I) and bigger 
if $|{\rm tr}e|\ll|{\rm dev}\,e|$ (i.e.\ Mode II), or very large if 
${\rm tr}e\ll-|{\rm dev}\,e|$ (compression leading to no fracture). 
Moreover, in the spirit of FFM from Remark~\ref{rem-stress-vs-engr}, one can 
consider energy in a finite neighbourhood of a current point, here 
split into the spherical and the shear parts to make $\zeta$ mode sensitive.
Of course, combination of both alternatives (i.e.\ also from 
Remark~\ref{rem-mode-sensitive}) is possible, too. On top of it,
one can also consider a combination with other dissipative processes
triggered only in Mode II, a prominent example being isochoric 
plasticity with hardening, cf.~Sect.\,\ref{sect-plast}.
Altogether, there are many parameters with clear physical interpretation
in the model to fit the model with many possible experiments in concrete
situations.
\end{remark}

\section{Various time discretisations}\label{sect-damage-discretisation}
In principle, the 2nd-order time derivative $\DDT u$ can be discretised by
2nd-order time differences, either as an explicit
\COL{(as in \eq{leapfrog} below)} or an implicit scheme.
This typically requires a fixed time step, and in the implicit variant
exhibits unacceptably
spurious numerical dissipation. Therefore, we 
avoid such discretisation here and work rather with the 1st-order system
\eqref{8-GSM1**-u-special} so that variable time-step is easily possible. 
Anyhow, for notational simplicity, we consider an equidistant partition of the
time interval $I=[0,T]$ with a fixed time step $\tau>0$ with $T/\tau\in\N$.

Considering some approximate values $\{u_\tau^k\}_{k=0,...,K}$ of the displacement $u$
with $K=T/\tau$, we define the piecewise-constant and the piecewise affine 
interpolants respectively by
\begin{subequations}\label{def-of-interpolants}
\begin{align}\label{def-of-interpolants-}
&&&
\overline{u}_\tau(t)= u_\tau^k,\qquad\ \
\underline u_\tau(t)= u_\tau^{k-1},\qquad\ \
\underline{\overline u}_\tau(t)=\frac12u_\tau^k+\frac12u_\tau^{k-1},
&&\text{and}
&&
\\&&&\label{def-of-interpolants+}
u_\tau(t)=\frac{t-(k{-}1)\tau}\tau u_\tau^k
+\frac{k\tau-t}\tau u_\tau^{k-1}
&&\hspace*{-7em}\text{for }(k{-}1)\tau<t\le k\tau.
\end{align}\end{subequations}
Similar meaning will have also $v_\tau$, $\overline v_\tau$, etc. 

\subsection{Implicit ``monolithic'' discretisation in time}\label{sec-monolit}

Some applications need to reflect the coupled character of the problem
in the truly coupled discrete fully-implicit scheme, in contrast to
the decoupled scheme considered in Sect.\,\ref{sec-stag} below.
This is indeed often solved in engineering, but only an approximate
solution can be expected by some iterative procedures.\footnote{\NEW{{Some
models are even formulated only in quasistatic time-discrete variants without
  having much chance to converge to some time-continuous problem; an example
  might be models with sharp interface between undamaged and partly damaged
  regions, as in \cite{AlJoGo09LSMN,XGNF17TDBF}.}}} Such schemes 
are known \COL{in engineering literature} under the adjective ``monolithic''
and the mentioned iterative solution is e.g.\ 
by the Newton-Raphson (or here equivalently SQP\,=\,sequential quadratic 
programming) method without any guaranteed convergence, however\NEW{, or
  alternating-minimization algorithm (AMA)\footnote{\COL{In the rate-independent
    quasistatic variant, AMA is similar the splitting scheme as in
    Sect.\,\ref{sec-stag} if the loading is modified as piecewise-constant
    in (rescaled) time except that the irreversibility constraint on the
 the damage profiles is up-dated differently.} \NEW{It was scrutinized e.g.\ in
      \cite{KneNeg17CAMS,MaMaPh16OMFG,Negr1QSEB} and used e.g.\ in
      \cite{LanRoy09VAFM}.}}}.
In general, such schemes even do not seem numerically stable because 
the a-priori estimates are not available. The semiconvexity here with 
respect to the $(H^1{\times}L^2$)-norm can be exploited 
provided the Kelvin-Voigt viscosity is used, 
as it is indeed considered in Sect.\,\ref{sect-damage-small-strains}
and \ref{sect-damage-PF-fracture}.

In contrast to the usual fully implicit scheme \COL{discretising
the inertial term by the second-difference formula
$\varrho(u_\tau^k{-}2u_\tau^{k-1}{+}u_\tau^{k-2})/\tau^2$ as e.g.\ in 
serving satisfactorily for analytical purpoces but causing an unacceptably
large spurious numerical dissipation, cf.\ e.g.\
\cite{BoLaRi11TDMD,LaOrSu10ESRM}}, we discretise the inertial 
part by the mid-point (Crank-Nicolson) formula rather than the backward Euler 
one in order to reduce unwanted numerical attenuation, and
we  use a semi-implicit 
(but not the fully implicit backward-Euler) 
formula for the visco-elastic stress
\COL{while} $\alpha$ is taken 
in an explicit way for the viscous part
in order to keep the 
variational structure of the incremental problems, cf.~\eq{potential} 
below, and to guarantee existence of the discrete solutions.
The resulted recursive coupled boundary-value problems here are:
\begin{subequations}\label{IBVP-damage-small-disc}
\begin{align}
&\!\!\!\frac{u_\tau^k{-}u_\tau^{k-1}}\tau=v_\tau^{k-1/2}\qquad\ \text{ with }\ \ 
v_\tau^{k-1/2}:=\frac{v_\tau^k{+}v_\tau^{k-1}}2,\\
&\!\!\!\varrho\frac{v_\tau^k{-}v_\tau^{k-1}}\tau-{\rm div}\,
\Big(
\bbD(\alpha_\tau^{k-1})e(v_\tau^{k-1/2})+\bbC(\alpha_\tau^{k})e(u_\tau^{k})
\Big)
=f_\tau^k
\label{IBVP-damage-disc-1}
\\&\qquad\qquad\qquad\qquad\qquad\text{ with }\ \ 
f_\tau^k:=\COL{\frac1\tau}\int_{(k-1)\tau}^{k\tau}\!\!f(t)\,\d t
,\ \text{ and }
\\
&
\!\!\!\partial\zeta
\Big(\frac{\alpha_\tau^k{-}\alpha_\tau^{k-1}}\tau\Big)
+
\frac12\bbC'(\alpha_\tau^{k})e(u_\tau^{k}){:}e(u_\tau^{k})
-{\rm div}(\kappa|\nabla\alpha_\tau^{k}|^{p-2}\nabla\alpha_\tau^{k})
\ni\fDAM'(\alpha_\tau^{k}) 
\label{IBVP-damage-disc-2}
\end{align}
\end{subequations}
considered on $\Omega$ while completed with the corresponding boundary 
conditions
\begin{subequations}\label{IBVP-damage-small-disc-BC}
\begin{align}\label{IBVP-damage-small-disc-BC-u}
&\Big(
\bbC(\alpha_\tau^{k})e(u_\tau^{k})
+\bbD(\alpha_\tau^{k-1})e(v_\tau^{k-1/2})
\Big)
\vec{n}=g_\tau^k\ \ \text{ and}
\\&\kappa\nabla\alpha_\tau^{k}{\cdot}\vec{n}=0,\ \
\text{ where }\ \ 
g_\tau^k:=\COL{\frac1\tau}\int_{(k-1)\tau}^{k\tau}\!\!g(t)\,\d t.
\end{align}
\end{subequations}
It is to be solved recursively for $k=1,...,T/\tau$, starting for $k=1$ with 
\begin{align}
u_\tau^0=u_0,\ \ \ \ \ \ 
v_\tau^0=v_0,
\ \ \ \ \ \ \alpha_\tau^0=\alpha_0.
\end{align}

In terms of the interpolants, see \eq{def-of-interpolants},
one can write the scheme \eqref{IBVP-damage-small-disc} more ``compactly'' as
\begin{subequations}\label{IBVP-damage-small-disc+-bis-}
\begin{align}
&\DT u_\tau=\underline{\overline v}_\tau\ \ \ \ \ \text{ and }\ \ \ \ \
\varrho\DT v_\tau-{\rm div}\,\big(\bbD(\underlinealpha_\tau)e(\underline{\overline v}_\tau)+\bbC(\overlinealpha_\tau)e({\overline u}_\tau)
\big)=\overline f_\tau\,,
\label{IBVP-damage-large-disc-1-bis-}
\\
&\pl\zeta\big(\DT\alpha_\tau\big)
+\frac12\bbC'(\overlinealpha_\tau)e(\overline u_\tau)\Colon e(\overline u_\tau)
-{\rm div}(\kappa|\nabla{\overlinealpha}_\tau|^{p-2}\nabla{\overlinealpha}_\tau)\ni
\fDAM'(\overlinealpha_\tau)\,.
\label{IBVP-damage-large-disc-2-bis-}
\end{align}
\end{subequations}
The boundary conditions \eq{IBVP-damage-small-disc-BC} can be written 
analogously.

Actually, we slightly modified the model used in 
Sections~\ref{sect-damage-small-strains} and 
Section~\ref{sect-damage-PF-fracture} by
considering a $p$-\COL{Laplacian}. For $p=2$, we obtain the previous
ansatz but for the convergence analysis we will need $p>d$.\footnote{See 
\eq{usage:p>d} below. In fact, the
presence of  $\bbD(\alpha_\tau^{k-1})$ instead of $\bbD(\alpha_\tau^k)$
brings difficulties in proving the strong convergence of rates, because
the 
analog of the argumentation 
used later in Sect.\,~\ref{sec-stag} does not work.
The mentioned non-quadratic modification of the gradient term is here
algorithmically tolerable because the strain energy 
$(e,\alpha)\mapsto\frac12\bbC(\alpha)e\Cdot e$ is not quadratic anyhow.} 
Using the ansatz 
\COL{\eq{eq6:AM-engr+}} with $\gamma$ smooth, positive, and strictly convex,
then $\frac12\bbC(\alpha)e{:}e+\frac12K|e|^2$ is convex for all $K$ large
enough. These underlying potentials
are strongly convex\footnote{To see it, one should analyze 
the Hessian on the the functional \eq{potential}, which is a bit
technical\COL{; cf.\ \cite{Roub??CTDD} for more details.}}
for the time-step $\tau>0$ small 
enough and, assuming also a conformal space discretisation, the iterative 
solvers have guaranteed convergence towards a unique (globally minimizing)
solution of the implicit scheme \eq{IBVP-damage-small-disc}. This
is satisfied 
for $\varphi$ from \eq{KV-damage-exa2} with $\eps=0$ or from 
\eq{eq6:AM-engr+}. The mentioned potential of the boundary-value
problem \eq{IBVP-damage-small-disc}--\eq{IBVP-damage-small-disc-BC} is 
\begin{align}\nonumber
&(u,\alpha)\mapsto\int_\Omega
\frac{\varrho}{2\tau}\Big|\frac{u{-}u_\tau^{k-1}}\tau{-}v_\tau^{k-1}\Big|^2\!
+\frac12\bbC(\alpha)e(u){:}e(u)-\fDAM(\alpha)+\tau\zeta\Big(\frac{\alpha{-}\alpha_\tau^{k-1}}\tau\Big)
\\[-.2em]&
\quad
+\frac{1}{2\tau}\bbD(\alpha_\tau^{k-1})e(u{-}u_\tau^{k-1}){:}e(u{-}u_\tau^{k-1})
+\frac\kappa p|\nabla\alpha(t)|^p
-f_\tau^k{\cdot}u\,\d x-\!\int_\Gamma\!g_\tau^k{\cdot}u\,\d S\,.
\label{potential}\end{align}
It is weakly lower semicontinuous on $H^1(\Omega;\R^d)\times H^1(\Omega)$
and coercive, so it serves also for proving existence of a weak solution 
to \eq{IBVP-damage-small-disc}--\eq{IBVP-damage-small-disc-BC}.
For any $(u_\tau^k,v_\tau^k,\alpha_\tau^k)\in H^1(\Omega;\R^d)\times 
L^2(\Omega;\R^d)\times W^{1,p}(\Omega)$ solving (in the usual weak sense)
the boundary value problem 
\eqref{IBVP-damage-small-disc}--\eqref{IBVP-damage-small-disc-BC},
the couple  $(u_\tau^k,\alpha_\tau^k)$ is a critical point of this functional. 
Also, conversely, any critical point $(u,\alpha)$ of \eqref{potential}
gives a weak solution  $(u_\tau^k,v_\tau^k,\alpha_\tau^k)$ 
to \eqref{IBVP-damage-small-disc}--\eqref{IBVP-damage-small-disc-BC} 
when putting $u_\tau^k=u$, $v_\tau^k=2(u_\tau^k{-}u_\tau^{k-1})/\COL{\tau}-v_\tau^{k-1}$,
and $\alpha_\tau^k=\alpha$.
For $\tau>0$ small enough, the mentioned convexity even ensures uniqueness
to this solution which is simultaneously a global minimizer of \eq{potential}.

The strategy \eqref{KV-visco-damageble-conv-strong} now uses 
the the piecewise affine and the piecewise constant
interpolants respectively as
\begin{align}
w_\tau=u_\tau+\taur v_\tau\ \ \ \text{ and }\ \ \ 
\overline w_\tau=\overline u_\tau+\taur\DT u_\tau=
\overline u_\tau+\taur\underline{\overline v}_\tau\,.
\end{align}
Then we can write the time-discrete approximation of the force equilibrium
\eq{IBVP-damage-large-disc-1-bis-}
as 
\begin{align*}
\frac{\varrho}{\taur}\DT w_\tau
-{\rm div}\big(\bbD_0e(\underline{\overline v}_\tau)
+\bbC(\overline\alpha_\tau)e(\overline w_\tau)\big)
=\overline f_\tau+\frac{\varrho}{\taur}\DT u_\tau+{\rm div}\big(
(\bbC(\overline\alpha_\tau)
{-}\bbC(\underline\alpha_\tau))e(\DT u_\tau)\big).
\end{align*}
We can test it by $\overline w_\tau$. By using in particular 
\begin{align}\nonumber
&\int_Q\!\DT w_\tau\!\cdot\!\overline w_\tau\,\d x\d t=
\int_Q\!(\DT u_\tau+\taur\DT v_\tau)
\!\cdot\!(\overline u_\tau+\taur\underline{\overline v}_\tau)\,\d x\d t
\\&\nonumber=
\int_Q\!(\DT u_\tau+\taur\DT v_\tau)
\!\cdot\!(\underline{\overline u}_\tau+\taur\underline{\overline v}_\tau)
+(\DT u_\tau+\taur\DT v_\tau)\!\cdot\!
(\overline u_\tau-\underline{\overline u}_\tau)\,\d x\d t
\\&=\int_\Omega\frac12|u_\tau(T)+\taur v_\tau(T)|^2
-\frac12|u_0+\taur v_0|^2\,\d x
+\!\!\!\!\lineunder{\frac\tau2\int_Q\!(\DT u_\tau{+}\taur\DT v_\tau)\!\cdot\!
\DT u_\tau\,\d x\d t}{$=\mathscr{O}(\tau)$}
\label{dwdt-vs-w}
\\[-3em]\nonumber\end{align}
we obtain the estimate
\begin{align}\nonumber
&\limsup_{\tau\to0}\int_Q\taur\bbD_0e(\underline{\overline v}_\tau)
{:}e(\underline{\overline v}_\tau)\,\d x\d t
\le\int_\Omega\frac\varrho{2\taur}|u_0{+}\taur v_0|^2+\frac12\bbD_0e(u_0){:}e(u_0)\,\d x
\\&\nonumber\
-\liminf_{\tau\to0}\bigg(\int_Q
\bbC(\overline\alpha_\tau)e(\overline w_\tau){:}e(\overline w_\tau)\,\d x\d t
+\int_\Omega\frac\varrho{2\taur}|u_\tau(T){+}\taur v_\tau(T)|^2+\frac12\bbD_0e(u_\tau(T)){:}e(u_\tau(T))\,\d x\bigg)
\\&\nonumber\ +\lim_{\tau\to0}\bigg(\int_Q\overline f_\tau{\cdot}\overline w_\tau
+(\bbC(\underline\alpha_\tau){-}\bbC(\overline\alpha_\tau)e(\DT u_\tau)
\!\colon\! e(\overline w_\tau)
\,\d x\d t
+\int_\Sigma\overline g_\tau{\cdot}\overline w_\tau\,\d S\d t
+\mathscr{O}(\tau)\bigg)
\\&\nonumber\ 
\le\int_\Omega\frac\varrho{2\taur}|u_0{+}\taur v_0|^2+\frac12\bbD_0e(u_0){:}e(u_0)\,\d x
+\int_Qf{\cdot}w-\bbC(\alpha)e(w){:}e(w)\,\d x\d t
\\&
\qquad\quad
+\int_\Omega\frac\varrho{2\taur}|u(T){+}\taur v(T)|^2+\frac12\bbD_0e(u(T)){:}e(u(T))\,\d x+\int_\Sigma g{\cdot}w\,\d S\d t
\ 
=\int_Q\taur\bbD_0e(v){:}e(v)\,\d x\d t,
\label{KV-visco-damageble-strong-Rothe}
\end{align}
where $\mathscr{O}(\tau)$ if from \eqref{dwdt-vs-w}.
The (last) equality in  \eqref{KV-visco-damageble-strong-Rothe} is 
due to the energy conservation in the limit equation 
\eqref{KV-visco-damageble-strong-}.
In \eqref{KV-visco-damageble-strong-Rothe}, we used also
\begin{align}\nonumber
&\bigg|\int_Q
\big(\bbC(\underline\alpha_\tau){-}\bbC(\overline\alpha_\tau)\big)
e(\DT u_\tau)\!\colon\! e(\overline w_\tau)\,\d x\d t\bigg|
\\&\ \ 
\le\big\|\bbC(\underline\alpha_\tau){-}\bbC(\overline\alpha_\tau)
\big\|_{L^\infty(Q;\R^{d^4})}
\big\|e(\DT u_\tau)\big\|_{L^2(Q;\R^{d\times d})}
\big\|e(\overline w_\tau)\big\|_{L^2(Q;\R^{d\times d})}\to 0
\label{usage:p>d}
\end{align}
Here we used 
the compact embedding of 
$L^\infty(I;W^{1,p}(\Omega))\,\cap\,H^1(I;L^2(\Omega))$ into $C(\overline Q)$ for 
$p>d$. This is actually one of the spot where $p>d$ is needed.

As we already know $e(\underline{\overline v}_\tau)\to e(v)$ weakly in 
$L^2(Q;\R^{d\times d})$, from \eqref{KV-visco-damageble-strong-Rothe}
we can see even the strong convergence. 
Since $e(\DT u_\tau)=e(\underline{\overline v}_\tau)$,
it also says $e(\DT u_\tau)\to e(\DT u)$ strongly, from which the desired strong
convergence $e(\overline u_\tau)\to e(u)$ needed for the limit passage in 
the damage flow rule follows.

\subsection{Fractional-step (staggered) discretisation}\label{sec-stag}

The damage problem typically involves the  \COL{stored} 
energies $\varphi=\varphi(e,\alpha)$ which are separately convex (or even 
separately quadratic). This encourages for an illustration of the 
\emph{fractional-step method},\index{fractional-step method!for damage}
also called \emph{staggered scheme}.\index{staggered scheme!for damage}  
In addition, to suppress a unwanted numerical attenuation within vibration, 
the time discretisation of the inertial term by the 
Crank-Nicholson scheme 
can also be considered, leading
to an energy-conserving discrete scheme. 

This falls into a broader class of the so-called HHT numerical integration 
methods devised by Hilber, Hughes, and Taylor \cite{HiHuTa77INDT},
generalizing the class of Newmark's methods \cite{Newm59MCSD},
widely used in engineering and computational physics. In fact, 
for a special choice of parameters,\footnote{In the 
standard notation used for the HHT-formula which uses three parameters, 
this special choice is $\alpha=\beta=1/2$ and $\gamma=1$.} the latter 
method gives the classical Crank-Nicolson scheme \cite{CraNic47PMNE} here 
applied to a transformed system of three 1st-order equations/inclusions 
\eq{8-GSM1*-special}. Actually, the Crank-Nicolson scheme
was originally devised for heat equation and later used for 
2nd-order problems in the form \eqref{8-GSM1**}, see 
e.g.\ \cite[Ch.6, Sect.9]{GlLiTr81NAVI}. It is different if applied to 
the dynamical equations transformed into the form \eq{8-GSM1*-special}; then
it is sometimes called just a central-difference scheme or 
generalized midpoint scheme, cf.\ e.g.\ \cite[Sect.\,12.2]{Wang07FSW} 
or \cite[Sect.\,1.6]{SimHug98CI}, respectively. For usage of Nemark's method
in dynamical damage see e.g.\ \cite{BVSH12,HofMie12CPFM,LMGP16NIDB,SWKM14PFAD}.

To allow \NEW{for} damage acting nonlinearly, we assume 
$\bbC(\cdot)$ and $\fDAM(\cdot)$ smooth and introduce the notation
\begin{align*}
\bbC_{ijkl}^\circ(\alpha,\tilde\alpha)
=\begin{cases}\displaystyle{\frac{\bbC_{ijkl}(\alpha)-\bbC_{ijkl}(\tilde\alpha)}{\alpha-\tilde\alpha}},&
\\
\bbC_{ijkl}'(\alpha)=\bbC_{ijkl}'(\tilde\alpha),&
\end{cases}
\fDAM^\circ(\alpha,\tilde\alpha)
=\begin{cases}\displaystyle{\frac{\fDAM(\alpha)-\fDAM(\tilde\alpha)}{\alpha-\tilde\alpha}}&\text{if }\alpha\ne\tilde\alpha,\\[.3em]
\fDAM'(\alpha)=\fDAM'(\tilde\alpha)&\text{if }\alpha=\tilde\alpha\,,
\end{cases}
\end{align*}
cf.\ e.g.\ \cite{CoMeSu11SAPF,RouPan17ECTD}.
Let us note that $\bbC^\circ(\alpha_\tau^k,\alpha_\tau^{k-1})
=\bbC'$ or $\fDAM^\circ(\alpha_\tau^k,\alpha_\tau^{k-1})
=\fDAM'$ if $\bbC(\cdot)$ or $\fDAM(\cdot)$ are affine.
It leads to the recursive boundary-value decoupled problems:
\begin{subequations}\label{IBVP-damage-small-disc+}
\begin{align}
&\frac{u_\tau^k{-}u_\tau^{k-1}}\tau=v_\tau^{k-1/2}:=\frac{v_\tau^k{+}v_\tau^{k-1}}2,
\\
&\varrho\frac{v_\tau^k{-}v_\tau^{k-1}}\tau
-{\rm div}\,\big(\bbD(\alpha_\tau^{k-1})
e(v_\tau^{k-1/2})+
\bbC(\alpha_\tau^{k-1})e(u_\tau^{k-1/2})\big)
=f_\tau^k\,,
\label{IBVP-damage-large-disc-1}
\\
&
\pl\zeta
\Big(\frac{\alpha_\tau^k-\alpha_\tau^{k-1}}\tau\Big)
+
\frac12\bbC^\circ(\alpha_\tau^k,\alpha_\tau^{k-1})e(u_\tau^k)\Colon e(u_\tau^k)
-{\rm div}(\kappa\nabla\alpha_\tau^{k-1/2})\ni\fDAM^\circ(\alpha_\tau^k,\alpha_\tau^{k-1})
\label{IBVP-damage-large-disc-2}
\end{align}
\end{subequations}
with $u_\tau^{k-1/2}:=\frac12u_\tau^k+\frac12u_\tau^{k-1}$ and
$\alpha_\tau^{k-1/2}:=\frac12\alpha_\tau^k+\frac12\alpha_\tau^{k-1}$,
considered on $\Omega$ while completed with the corresponding boundary 
conditions discretized analogously. It is to be solved recursively for 
$k=1,...,T/\tau$, starting with 
\begin{align}\label{IBVP-damage-small-disc-IC}
u_\tau^0=u_0,\ \ \ \ \ \ 
v_\tau^0=v_0,\ \ \ \ \ \ \alpha_\tau^0=\alpha_0,
\end{align}
and solving alternately (\ref{IBVP-damage-small-disc+}a,b) and 
\eq{IBVP-damage-large-disc-2}. Both these boundary-value problems have their 
own potentials.

In terms of the interpolants, see \eq{def-of-interpolants},
one can write the scheme \eqref{IBVP-damage-small-disc+} more ``compactly'' as
\begin{subequations}\label{IBVP-damage-small-disc+-bis}
\begin{align}
&\DT u_\tau=\underline{\overline v}_\tau\ \ \ \ \ \text{ and }\ \ \ \ \
\varrho\DT v_\tau-{\rm div}\,\big(\bbD(\underlinealpha_\tau)
e(\underline{\overline v}_\tau)+\bbC(\underlinealpha_\tau)e(\underline{\overline u}_\tau)
\big)=\overline f_\tau\,,
\label{IBVP-damage-large-disc-1-bis}
\\
&\pl\zeta\big(\DT\alpha_\tau\big)
+\frac12\bbC^\circ(\overlinealpha_\tau,\underlinealpha_\tau)e(\overline u_\tau)\Colon e(\overline u_\tau)
-{\rm div}(\kappa\nabla\underline{\overlinealpha}_\tau)\ni
\fDAM^\circ(\overlinealpha_\tau,\underlinealpha_\tau)\,.
\label{IBVP-damage-large-disc-2-bis}
\end{align}
\end{subequations}
The boundary conditions can be written analogously.
The basic energetic test of \eq{IBVP-damage-large-disc-1-bis} is to be done by 
$\DT u_\tau=\underline{\overline v}_\tau$ and of \eq{IBVP-damage-large-disc-2-bis}
by $\DT\alpha_\tau$. We can use a binomial formula several times,
in particular for 
\begin{subequations}\begin{align}
&\varrho\frac{v_\tau^k{-}v_\tau^{k-1}}\tau\Cdot\frac{v_\tau^k{+}v_\tau^{k-1}}2
=\frac{\frac12\varrho|v_\tau^k|^2{-}\frac12\varrho|v_\tau^{k-1}|^2}\tau,
\\&\kappa\nabla\frac{\alpha_\tau^k+\alpha_\tau^{k-1}}2\Cdot\nabla
\frac{\alpha_\tau^k-\alpha_\tau^{k-1}}\tau
=\frac{\frac12\kappa|\nabla\alpha_\tau^k|^2
-\frac12\kappa|\nabla\alpha_\tau^{k-1}|^2}\tau,\ \ \text{ and}
\\&\nonumber
\bbC(\alpha_\tau^{k-1})\!\!\!\!\!\!\!\lineunder{e(u_\tau^{k-1/2})\Colon e(v_\tau^{k-1/2})}{$=
\frac{e(u_\tau^{k})\Colon e(u_\tau^{k})-e(u_\tau^{k-1})\Colon e(u_\tau^{k-1})}{2\tau}$}\!\!\!\!\!\!
+\frac12\!\lineunder{\frac{\alpha_\tau^k-\alpha_\tau^{k-1}}\tau
\bbC^\circ(\alpha_\tau^k,\alpha_\tau^{k-1})}{$\ \ \ =(\bbC(\alpha_\tau^k)-\bbC(\alpha_\tau^{k-1}))/\tau$}
e(u_\tau^k)\Colon e(u_\tau^k)
\\&\qquad\qquad=\Big(\frac12\bbC(\alpha_\tau^k)e(u_\tau^{k})\Colon e(u_\tau^{k})
-\frac12\bbC(\alpha_\tau^{k-1})e(u_\tau^{k-1})\Colon e(u_\tau^{k-1})\Big)/\tau\,;
\end{align}\end{subequations} 
note that we have enjoyed 
the cancellation of the terms 
$\pm\frac12\bbC(\alpha_\tau^{k-1})e(u_\tau^k)\Colon e(u_\tau^k)$, cf.\ also \cite{Roub17ECTD}.
Thus we obtain the discrete analog of energy equality \eq{KV-damage-engr}:
\begin{align}\nonumber
&\int_\Omega\frac\varrho2|\DT u_\tau(t)|^2+
\frac12\bbC(\alpha_\tau(t))e(u_\tau(t))\Colon e(u_\tau(t))-\fDAM(\alpha(t))
+\frac\kappa 2|\nabla\alpha(t)|^2\,\d x
+\int_0^t\!\!\int_\Omega\bbD(\underlinealpha_\tau) e(\DT u_\tau)\Colon e(\DT u_\tau)
+\DT\alpha_\tau\pl\zeta(\DT\alpha_\tau)\,\d x\d t
\\[-.4em]&
=\int_\Omega\frac\varrho2|v_0|^2+\frac12\bbC(\alpha_0)e(u_0)\Colon e(u_0)
-\fDAM(\alpha_0)
+\frac\kappa 2|\nabla\alpha_0|^2\,\d x
+\int_0^t\!\!\int_\Omega\overline f_\tau\Cdot\DT u_\tau\,\d x\d t+
\int_0^t\!\!\int_\Gamma \overline g_\tau\Cdot\DT u_\tau\,\d\NEW{S}\d t
\label{KV-damage-engr-disc}
\end{align}
at each mesh point $t=k\tau$ with $k\in\{0,...,T/\tau\}$. Let us note that this
is indeed an equality, not only an estimate. This discrete energy conservation
can advantageously be used to check a-posteriori correctness of a 
computational code.

We introduce the variables $w_\tau^k=u_\tau^k+\taur v_\tau^k$ for 
$k\in\{0,...,T/\tau\}$ and the corresponding interpolants 
$w_\tau=u_\tau+\taur v_\tau$ and $\underline{\overline w}_\tau
=\underline{\overline u}_\tau+\taur\underline{\overline v}_\tau$.
Likewise \eq{KV-visco-damageble-strong-}, we can rewrite 
\eq{IBVP-damage-large-disc-1-bis} as 
\begin{align}
&
\frac{\varrho}{\taur}\DT w_\tau-{\rm div}\big(\bbD_0e(\DT u_\tau)+
\bbC(\underlinealpha_\tau)e(\underline{\overline w}_\tau)\big)
=\overline f_\tau+\frac{\varrho}{\taur}\DT u_\tau.
\label{semiimplicit-w}\end{align}
To replicate the strategy \eq{KV-visco-damageble-conv-strong}, we use 
a test of \eq{semiimplicit-w} by $\underline{\overline w}_\tau-w$
and the calculus
\begin{align}\nonumber
\int_Q\varrho(\DT w_\tau-\DT w)\Cdot(\underline{\overline w}_\tau-w)\,\d x\d t
  &
=
\int_\Omega\frac\varrho2|w_\tau(T)-w(T)|^2\,\d x
+\int_Q\varrho(\DT w_\tau-\DT w)\Cdot(\underline{\overline w}_\tau-w_\tau)\,\d x\d t
\\[-.3em]&
=
\int_\Omega\frac\varrho2|w_\tau(T)-w(T)|^2\,\d x
-\int_Q\!\!\!\!\!\!\!\lineunder{\,\varrho\DT w\Cdot(\underline{\overline w}_\tau-w_\tau)}{$\to0$ in $L^1(Q)$ weakly}\!\!\!\!\!\!\d x\d t
\end{align}
because $\int_0^T\DT w_\tau\Cdot(\underline{\overline w}_\tau-w_\tau)\,\d t=0$
a.e.\ on $\Omega$.
Similarly, still we use the calculus
\begin{align}\nonumber
&\int_Q\bbD_0 e(\underline{\overline v}_\tau{-}v)\Colon e(\underline{\overline w}_\tau{-}w)\,\d x\d t=\int_Q\bbD_0 e(\DT u_\tau{-}\DT u)\Colon e(\underline{\overline u}_\tau{-}u)
+\taur\bbD_0 e(\underline{\overline v}_\tau{-}v)\Colon e(\underline{\overline v}_\tau{-}v)
\,\d x\d t
\\[-.3em]&\qquad
=\int_\Omega\bbD_0 e(u_\tau(T){-}u(T))\Colon e(u_\tau(T){-}u(T))\,\d x
+\int_Q\taur\bbD_0 e(\underline{\overline v}_\tau{-}v)\Colon e(\underline{\overline v}_\tau{-}v)-\!\!\!\lineunder{\bbD_0 e(\DT u)\Colon e(\underline{\overline u}_\tau{-}u)}{$\to0$ in $L^1(Q)$ weakly}\!\!\!
\,\d x\d t\,.
\end{align}
Thus we obtain the strong convergence $e(\underline{\overline v}_\tau)
=e(\DT u_\tau)\to e(\DT u)$ in $L^2(Q;\Rsym)$, and thus also 
$e(\overline u_\tau)\to e(u)$ needed to pass to the limit in
\eq{IBVP-damage-large-disc-2-bis}.

When using the separately quadratic ansatz \eq{eq6:AM-engr+} and when
combined the time discretisation \eq{IBVP-damage-small-disc+}
with P1 finite-element space discretisation, it gives an alternating 
linear-quadratic programming problems and thus very efficient
numerical algorithms; in fact, it can be implemented without any 
iterative procedure needed, and the energy balance \eq{KV-damage-engr-disc}
is satisfied exactly up to only round-off errors. 

\begin{figure}
\begin{center}
  \centerline{\includegraphics[height=1.5\textwidth]{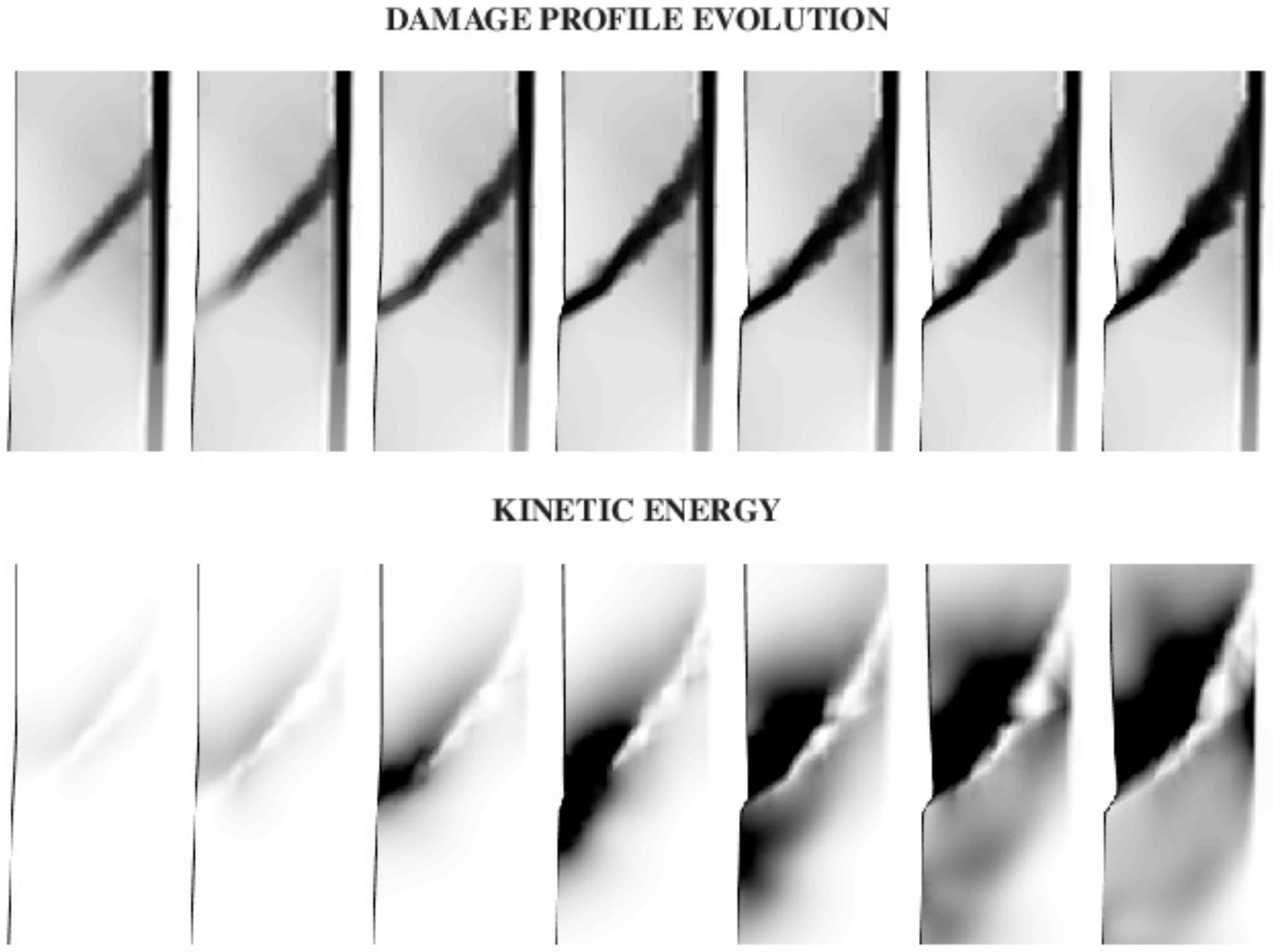}}
\end{center}
\vspace*{-40em}
\caption{\sl Simulations of a rupture in a two-dimensional specimen loaded 
by tension in a vertical direction, modelled by the phase-field crack 
approximation, and subsequent emission of an elastic wave. Seven selected 
snapshots are depicted. The decoupled energy-preserving time discretisation 
and P1-finite elements have been used.
\newline\ \hspace*{19em}Courtesy of Roman Vodi\v cka 
(Technical University Ko\v sice, Slovakia)}
\label{fig-simul}
\end{figure}

Let us briefly illustrate this algorithm on a 2-dimensional computational 
experiment considering an isotropic material occupying a 
rectangular domain $\Omega$.
The left side of this rectangular vertically stretched specimen
is left free while the right-hand side is allow\NEW{ed} to slide. 
This asymmetry also causes a slight asymmetry of the solution and not 
completely straight fracture line, cf.\ Figure~\ref{fig-simul}. 
Although the discretisation scheme is unconditionally convergent, to
see reasonable numerical results,
one should respect the maximal wave speed by choosing reasonably small
time step, cf. the CFL-condition in the following Sect.~\ref{sec-explicit}.
For details 
about the implementation and data and more complete presentation of an overall 
experiments we refer to \cite{RouVod??MMSS}.

\subsection{
  Explicit time discretisation \COL{outlined}}\label{sec-explicit}

Implicit schemes from Sect.\,\ref{sec-monolit} and \ref{sec-stag}
are not causal and not much efficient for real wave propagation calculations
\COL{usually containing higher frequencies 
in comparison with mere vibrations}. For this, more often, explicit schemes are 
used for real wave calculations, at least in linear
\COL{elastodynamic} models. These time-discretisation schemes work only if
combined with space discretisation.

Disregarding damage and the viscous rheology, one efficient option often 
considered for waves in purely elastic materials is a so-called leapfrog 
scheme (also known as Verlet's integration), i.e.\ central differences for
the kinetic term
\begin{align}\label{leapfrog}
\mathscr{T}'\frac{v_\tau^{k+1/2}-v_\tau^{k-1/2}}{\tau}
+\pl^{}_u\mathscr{E}_h(u_\tau^{k})=\mathscr{F}_\tau^{k}\ \ \ \ \text{ with }\ \
v_\tau^{k+1/2}=\frac{u_\tau^{k+1}{-}u_\tau^k}\tau.
\end{align}
The test by $v_\tau^{k+1/2}$ leads to a slightly twisted energy 
(im)balance:
$$
\mathcal{T}\big(v_\tau^{k+1/2}\big)
+\frac12\big\langle\pl^{}_u\mathscr{E}_h(u_\tau^{k}),u_\tau^{k+1}\big\rangle
=\mathcal{T}\big(v_\tau^{k-1/2}\big)
+\frac12\big\langle\pl^{}_u\mathscr{E}_h(u_\tau^{k-1}),u_\tau^k\big\rangle
+\big\langle\mathscr{F}_\tau^{k},v_\tau^{k+1/2}\big\rangle.
$$
This gives a correct kinetic energy but the stored energy is 
correct only asymptotically under the Courant-Friedrichs-Lewy
(so-called CFL) condition \cite{CoFrLe28PDMP}, which needs also 
a space discretisation (here indicated by the abstract ``mesh parameter'' 
$h>0$) and the time step sufficiently small with respect to $h$; typically
$\tau<Vh$ with $V$ the maximal speed of arising waves
if $h$ has the meaning of a size of the largest element in a finite-element 
discretisation.

\COL{The option \eq{leapfrog} does not seem directly amenable for being merged
  with the damage evolution.} Another option
\COL{relies on the reformulation of the elastodynamics in terms
  of velocity and stress, i.e.\ in terms of $v=\DT u$ and of
  the stress $\sigma:=
  \bbC e(u)$, eliminating the displacement $u$.
  We thus have in mind the system
\begin{subequations}\label{IBVP+}\begin{align}\label{IBVP1+}
    &&&  \DT\sigma=\bbC e(v)\ \ \ \text{ and }\ \ \
    \varrho\DT v-{\rm div}\,\sigma
    =f
    &&\text{in }\ Q,&&
  \\&&&\label{IBVP2+}
  \DT\sigma\vec{n} 
  =\DT g&&\text{on }\ \varSigma,
  \\&&&\label{IBVP3+}
  v|_{t=0}=v_0,\ \ \ \sigma|_{t=0}=\sigma_0:=\bbC e(u_0)&&\text{\NEW{in} }\ \varOmega.
\end{align}\end{subequations}
The explicit staggered (called also ``leap-frog'') time-discretisation can
now be done as
\begin{subequations}\label{IBVP+disc}\begin{align}\nonumber\\[-2.5em]
    \label{IBVP1+disc}
    &&& 
    \frac{\sigma_\tau^{k}-\sigma_\tau^{k-1}}\tau=\bbC e(v_\tau^{k-1})
    \ \ \text{ and }\ \ 
    \varrho\frac{
      v_\tau^{k}-v_\tau^{k-1}}\tau
    -{\rm div}
    \sigma_\tau^{k}
    =f_\tau^{k}
    &&\text{in }\ Q,&&
  \\[-.5em]&&&\label{IBVP2+disc}
\frac{\sigma_\tau^{k}-\sigma_\tau^{k-1}}\tau\vec{n}
  =\frac{g_\tau^{k}-g_\tau^{k-1}}\tau&&\text{on }\ \varSigma,
  \\&&&\label{IBVP3+disc}
  v_\tau^0=v_0,\ \ \ \sigma_\tau^0=\sigma_0:=\bbC e(u_0)&&\text{\NEW{in} }\ \varOmega.
\end{align}\end{subequations}
Let us note that \eq{IBVP1+disc} is decoupled, i.e.\
one is first \NEW{to compute}
$\sigma_\tau^{k}$ 
and then $v_\tau^{k}$.
Averaging the second equation in \eq{IBVP1+disc} at level $k$ and $k{-}1$
and testing it $v_\tau^{k-1}$ while  
testing the first  equation in \eq{IBVP1+disc} by
$(\sigma_{\tau h}^{k}+\sigma_{\tau h}^{k-1})/2$, we obtain the approximate
energy balance as
\begin{align}
&\frac12\big\langle\mathcal{T}'v_\tau^{k},v_\tau^{k-1}\big\rangle
+\varPhi_h(\sigma_\tau^{k})
\label{...}
=\frac12\big\langle\mathcal{T}'v_\tau^{k-1},v_\tau^{k-2}\big\rangle
+\varPhi_h(\sigma_\tau^{k-1})
+\big\langle\mathscr{F}_\tau^{k},v_\tau^{k-1}\big\rangle\,,
\end{align}
with $\varPhi$ the stored energy expressed in terms of stress.
Now the stored energy is correct while the kinetic energy needs the 
CFL-condition, cf.\ \cite{ByJoTs02NFMF}. In contrast with \eq{leapfrog},
this option is more compatible with possible 
enhancement of the stored energy by internal parameters as e.g.\
damage. 

Assuming $\bbC(\alpha)=\gamma(\alpha)\bbC_1$ as in \eq{eq6:AM-engr+},
we consider the energy $\varPhi=\varPhi(\varsigma,\alpha)$
with a ``proto-stress'' $\varsigma=\bbC_1 e(u)$ and with
$\varPhi(\varsigma,\alpha)
=\int_\varOmega\frac12\gamma(\alpha)\bbC_1^{-1}\varsigma{:}\varsigma
-\fDAM(\alpha)+\frac\kappa2|\nabla\alpha|^2\,\d x$; for a general concept
see \cite{RoPaTs??ETDE} although, in damage mechanics, this proto-stress
is also called an effective stress, having a specific mechanical meaning
\cite{Rabo69CPSM}. An important trick is that the proto-stress does not
explicitly involve $\alpha$ and its time derivate does not lead to $\DT\alpha$.
The system \eq{IBVP+} enhnaced by damage like \eq{8-GSM1**-z-special}
then looks as
\begin{subequations}\label{IBVP+dam}\begin{align}\label{IBVP1++dam}
    &&& \DT\varsigma=\bbC_1 e(v)\ \ \ \text{ and }\ \ \
    \varrho\DT v-{\rm div}\,\sigma
    =f\ \ \ \text{ with }\ \ \
\sigma=\gamma(\alpha)\varsigma
    &&\text{in }\ Q,&&
    \\\label{IBVP2++dam++}
    &&&\pl\zeta(\DT\alpha)+
    \frac12
    \gamma'(\alpha)\bbC_1^{-1}\varsigma{:}\varsigma
-{\rm div}\big(\kappa
\nabla \alpha\big)
\ni\fDAMprime(\alpha)
&&\text{in }\ Q,
  \\&&&\label{IBVP2++dam}
  \DT\sigma\vec{n}
  =\DT g\ \ \text{ and }\ \ \kappa
\nabla \alpha\cdot\vec{n}=0&&\text{on }\ \varSigma,
  \\&&&\label{IBVP3+dam}
  v|_{t=0}=v_0,\ \ \ \sigma|_{t=0}=\sigma_0:=\bbC e(u_0),
  \ \ \ \alpha|_{t=0}=\alpha_0&&\text{\NEW{in} }\ \varOmega.
\end{align}\end{subequations}
Applying the staggered discretisation like in Sect.~\ref{sec-stag},
we
obtain a 3-step scheme:
\begin{subequations}\label{IBVP+disc+}\begin{align}\label{IBVP1+disc+}
    &&& 
    \frac{\varsigma_\tau^{k}-\varsigma_\tau^{k-1}}\tau=\bbC_1 e(v_\tau^{k-1})&&\text{in }\ Q\,,&&
    \\
    &&&\pl\zeta
\Big(\frac{\alpha_\tau^k-\alpha_\tau^{k-1}}\tau\Big)
+
\frac12\gamma^\circ(\alpha_\tau^k,\alpha_\tau^{k-1})\bbC_1^{-1}\varsigma_\tau^k
\Colon\varsigma_\tau^k
-{\rm div}(\kappa\nabla\alpha_\tau^{k-1/2})\ni\fDAM^\circ(\alpha_\tau^k,\alpha_\tau^{k-1})&&\text{in }\ Q\,,&&
\label{IBVP3++disc++}
\\&&&
   \varrho\frac{
      v_\tau^{k}-v_\tau^{k-1}}\tau
    -{\rm div}
    \sigma_\tau^{k}
    =f_\tau^{k}\ \ \ \text{ with }\ \ \
    \sigma_\tau^{k}=\gamma(\alpha_\tau^k)\varsigma_\tau^{k}&&\text{in }\ Q\,,&&
\end{align}\end{subequations}
to be completed by the respective boundary conditions. 

The analysis of \eq{IBVP+disc+} is however rather nontrivial and the analog
of \eq{...} with the corresponding damage terms like in \eq{KV-damage-engr-disc}
contains still some other term vanishing in the limit under the
CFL condition, cf.\ \cite{RoPaTs??ETDE} for details. Even more, as
there is no Kelvin-Voigt viscosity which would be troublesome for
such explicit discretisation, one needs still some higher-order
gradient term not subject to damage and acting on $\varsigma$
to guarantee convergence of such a scheme; cf.\ also
\cite[Sect.7.5.3]{KruRou19MMCM}.

To 
\NEW{conclude}, it should be mentioned that a really efficient
(i.e.\ explicit) numerical scheme with granted stability and convergence
for the simple inviscid or viscous material undergoing damage does not seem to 
be devised so far.

}

\section{Concluding remarks -- some modifications}
\label{sect-damage-generalization}
Many other phenomena can be combined with the plain damage in
the Kelvin-Voigt vicoelastic model considered so far. 
Typically \NEW{one can} think about more complicated viscoelastic rheologies,
involving possibly some inelastic processes as plasticity,
which will be in a simple variant in Sect.~\ref{sect-plast}. 

Also, some diffusant (like water in poroelastic rocks or hydrogen in metals or 
some solvent in polymers) can propagate through the bulk by a Fick/Darcy law, 
interacting with mechanical properties including fracture toughness.
Of course, full thermodynamical context should involve heat production and  
transfer through the Fourier law. Here we only refer to \cite{Roub17ECTD} 
where a staggered energy-conserving time discretisation like in 
Sect.~\ref{sec-stag} is devised. Damage with plasticity accompanied by heat 
production and heat transfer allows for fitting to the popular 
rate-and-state-dependent friction model \cite{Roub14SNRSD}.

Moreover, the plain models from 
Sections~\ref{sect-damage-small-strains}--\ref{sect-damage-PF-fracture} 
together with all these extensions can be considered \COL{within} the large strains, 
too. We will 
\COL{outline} it Sect.~\ref{sect-large}.

\subsection{Combination with creep or plasticity}\label{sect-plast}
\index{plasticity}

The Kelvin-Voigt rheology is mathematically the most basic viscoelastic 
rheology of parabolic type. In particular from the wave-propagation
viewpoint, physically more natural is that Maxwell rheology but it is 
rather hyperbolic and mathematically troublesome if accompanied with 
inelastic processes like damage. A certain reasonable compromise it the
Jeffreys' rheology combining the Norton-Hoff (also called Stokes) 
and Kelvin-Voigt rheology in series (or alternatively Maxwell's and 
Norton-Hoff's rheology in parallel). It can capture creep effects,
which have sense in the shear part rather than the spherical part.

\begin{figure}
\begin{center}
      {\includegraphics[height=1.6\textwidth]{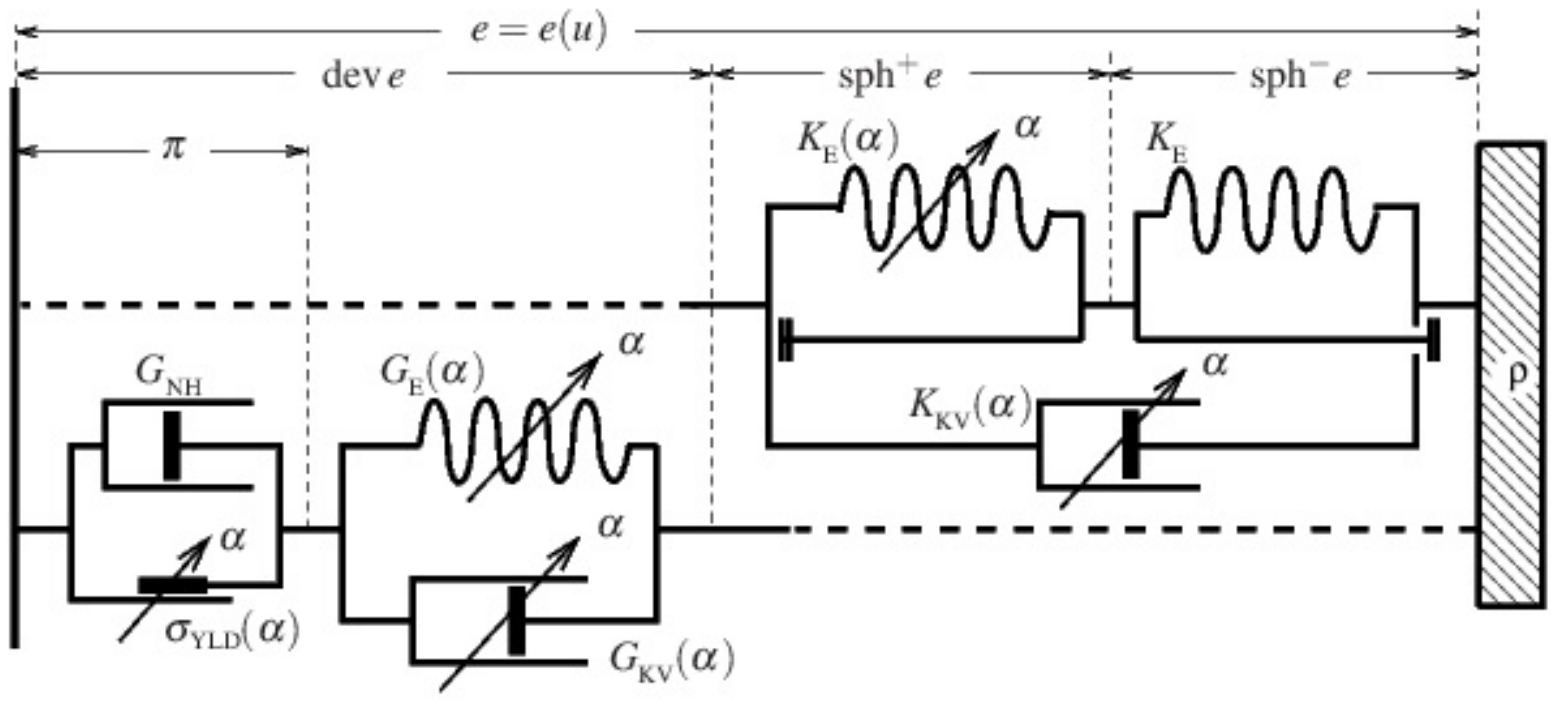}}
\end{center}
\vspace*{-55em}
\caption{\sl Schematic diagram for the viscoelastic Jeffreys rheology
(if $\SYLD=0$) 
which is subjected to damage $\alpha$ in the deviatoric part except 
undamageable creep (the $G_{_{\rm NH}}$-dashpot) while the Kelvin-Voigt rheology 
in the spherical (volumetric) part is subjected to damage only under 
tension but not compression. For $\SYLD>0$, it models (visco)plasticity.
Evolution of damage is not depicted.}
\label{fig-rheology}
\end{figure}

Instead (or in addition) to the linear Norton-Hoff dumper in the 
shear part, one can consider also the activated plastic element.
The schematic rheological model is depicted in Figure~\ref{fig-rheology},
\COL{distinguishing} also the compression and the tension in the  spherical part.


The additional dissipation due to isochoric plastification is then achieved 
when damage is performed in a shear mode (i.e.\ Mode II) comparing to damage 
by opening (i.e.\ Mode I) where plastification is not triggered. When 
considering the isotropic stored energy \eq{KV-damage-exa2} with
damage \NEW{without any hardening-like effects} (i.e.\ linearly-depending $K_{_{\rm E}}(\alpha)=\alpha K$ 
and $G_{_{\rm E}}(\alpha)=\alpha G$) and with 
$\fDAM(\alpha)$ also linear and the elastic strain $e-\pi$ in 
place of the total strain $e$ combined with the isotropic hardening 
with 
$\SEFS{}=0$, we altogether arrive at the model governed by 
\begin{subequations}\label{plast-dam-syst-mixed}
\begin{align}
\nonumber
&\varphi_\text{\sc e}(e,\alpha,\pi,\nabla\alpha,\nabla\pi)=
\frac d2\Big(K_{_{\rm E}}|{\rm sph}^-e|^2
+K_{_{\rm E}}(\alpha)
|{\rm sph}^+e|^2\Big)
+G_{_{\rm E}}(\alpha)|{\rm dev}e{-}\pi|^2
\\&\qquad\qquad\qquad\qquad\qquad\quad
-\fDAMprime(\alpha)
+\frac12H|\pi|^2
+\frac{\kappa_1}2|\nabla\pi|^2
+\frac{\kappa_2}2|\nabla\alpha|^2,
\label{plast-dam-syst-mixed-chi}
\\&
\zeta(\alpha;\DT e,\DT\alpha,\DT\pi)=
\frac d2 K_\text{\sc kv}(\alpha)|{\rm sph}\,\DT e|^2+G_\text{\sc kv}
(\alpha)|{\rm dev}\DT e{-}\DT\pi|^2
+\frac12G_\text{\sc nh}|\DT\pi|^2
\label{plast-dam-syst-mixed-zeta}
+\SYLD(\alpha)|\DT\pi|+\delta_{[0,
+\infty)}^*(\DT\alpha)
+\frac{\nu}2\DT\alpha^2,
\end{align}\end{subequations}
where the specific dissipation potential now contains another damper 
$G_\text{\sc nh}$ which facilitates to the Jeffrey's model in the shear 
part and a yield stress $\SYLD{}\ge0$ possibly depending on damage, which 
can model activated inelastic plastic response. Starting from undamaged 
material, the energy needed (dissipated) by damaging in opening without 
plastification is just the toughness $\Frakg_{\rm c}:=\fDAM'(1)$, while 
in shearing mode it is larger, namely $\Frakg_{\rm c}+\SYLD(\sqrt{2G\Frakg_{\rm c}}
-\SYLD)/H$ provided the parameters are tuned in a way to satisfy 
$\sqrt{G\Frakg_{\rm c}/2}<\SYLD\le\sqrt{2G\Frakg_{\rm c}}$.
This was first devised for an interfacial delamination model 
\cite{RoMaPa13QMMD}, being inspired just by such bulk plasticity.

Let us illustrate the staggered scheme \eq{IBVP-damage-small-disc+-bis} 
in the case of a linearly responding material, i.e.\ 
$K_{_{\rm E}}|{\rm sph}^-e|^2
+K_{_{\rm E}}(\alpha)|{\rm sph}^+e|^2$ is simplied to 
$K_{_{\rm E}}(\alpha)|{\rm sph}\,e|^2$ in \eq{plast-dam-syst-mixed-chi},
denoting $\bbC_{ijkl}(\alpha)=K_\text{\sc e}(\alpha)\delta_{ij}\delta_{kl}
+G_\text{\sc e}(\alpha)(\delta_{ik}\delta_{jl}
+\delta_{il}\delta_{jk}-\frac2d\delta_{ij}\delta_{kl})$ 
and $\bbD_{ijkl}(\alpha)=K_\text{\sc kv}(\alpha)\delta_{ij}\delta_{kl}
+G_\text{\sc kv}(\alpha)(\delta_{ik}\delta_{jl}
+\delta_{il}\delta_{jk}-\frac2d\delta_{ij}\delta_{kl})$ 
with $\delta$ standing for the Kronecker symbol.
More specifically, introducing a notation for the elastic 
strain $e_{\rm el}=e(u){-}\pi$ and its discretisation 
 $e_{\rm el,\tau}=e(u_\tau){-}\pi_\tau$ and 
$\underline{\overline e}_{\rm el,\tau}
=e(\underline{\overline u}_\tau){-}\underline{\overline\pi}_\tau$,
the system \eq{IBVP-damage-small-disc+-bis} can be expanded as
\begin{subequations}\label{IBVP-damage-small-disc+plast}
\begin{align}
&\DT u_\tau=\underline{\overline v}_\tau\ \ \ \ \ \text{ and }\ \ \ \ \
\varrho\DT v_\tau-{\rm div}\big(\bbD(\underlinealpha_\tau)
\DT e_{\rm el,\tau}
+\bbC(\underlinealpha_\tau)
\underline{\overline e}_{\rm el,\tau}
\big)=\overline f_\tau\,,
\label{IBVP-damage-large-disc-1-plast}
\\
&
\SYLD(\underlinealpha_\tau){\rm Dir}(\DT\pi_\tau)+H\underline{\overline \pi}_\tau
-{\rm div}(\kappa_1\nabla\underline{\overline\pi}_\tau)
\ni{\rm dev}
\big(\bbD(\underlinealpha_\tau)\DT e_{\rm el,\tau}
+\bbC(\underlinealpha_\tau)\underline{\overline e}_{\rm el,\tau}
\big),
\label{IBVP-damage-large-disc-12-plast}
\\&
\pl\zeta\big(\DT\alpha_\tau\big)
+\frac12\bbC^\circ(\overlinealpha_\tau,\underlinealpha_\tau)
\overline e_{\rm el,\tau}\Colon\overline e_{\rm el,\tau}
-{\rm div}(\kappa_2\nabla\underline{\overlinealpha}_\tau)\ni
\fDAM^\circ(\overlinealpha_\tau,\underlinealpha_\tau)\,.
\label{IBVP-damage-large-disc-3-plast}
\end{align}
\end{subequations}
The boundary conditions can be written analogously. Now, the splitting
during the recursive time-stepping procedure concerns separately
(\ref{IBVP-damage-small-disc+plast}a,b) and 
\eq{IBVP-damage-large-disc-3-plast}, both these boundary-value problems 
at particular time levels having a potential. The basic energy estimates
can be obtained by testing the particular equations/inclusions in 
\eq{IBVP-damage-small-disc+plast} subsequently by $v_\tau$, $\DT\pi_\tau$,
and $\DT\alpha_\tau$, using the quadratic trick several times 
(e.g.\ for $\varrho\DT v_\tau\cdot\underline{\overline v}_\tau
=\frac{\partial}{\partial t}\frac12|\underline{\overline v}_\tau|^2$
a.e.\ o $Q$) and the cancellation of the terms $\pm\bbC(\underlinealpha_\tau)
\overline e_{\rm el,\tau}\Colon\overline e_{\rm el,\tau}$ arising by these tests.

For the strong convergence of $\overline e_{\rm el,\tau}$, instead of the 
strategies \eq{KV-visco-damageble-conv-strong} or 
\eq{KV-damage-strong-e(t)}, we now rely rather on the test 
of (\ref{IBVP-damage-small-disc+plast}a,b) respectively by 
$\DT u_\tau{-}\DT u$ and $\DT\pi_\tau{-}\DT\pi$.
We need $\bbC$ monotone (nondecreasing) with respect
to the L\"owner ordering, and we use the unidirectionality
of the damage evolution, i.e.\ $\DT\alpha_\tau\le0$.
We first approximate the limit $u$ and $\pi$,
defining $\widetilde u_\tau^k:=\COL{\frac1\tau}\int_{(k-1)\tau}^{k\tau}u(t)\,\d t$
and  $\widetilde\pi_\tau^k:=\COL{\frac1\tau}\int_{(k-1)\tau}^{k\tau}\pi(t)\,\d t$,
and then put
$\varepsilon_{\rm el,\tau}^k:=e(\widetilde u_\tau^k)-\widetilde\pi_\tau^k$.
Then also the interpolants $\varepsilon_{\rm el,\tau}$ and 
$\underline{\overline\varepsilon}_{\rm el,\tau}$ which both converges to 
$e_{\rm el}$ strongly. 
This approximation allows us to estimate\footnote{Here,
abbreviating $E_\tau^k=e_{\rm el,\tau}^k-\varepsilon_{\rm el,\tau}^k$, we used the 
algebra 
\begin{align*}
&\frac12\bbC(\alpha_\tau^k)E_\tau^k\Colon E_\tau^k
-\frac12\bbC(\alpha_\tau^{k-1})E_\tau^{k-1}\Colon E_\tau^{k-1}
\\&\qquad=\frac12(\bbC(\alpha_\tau^k)-\bbC(\alpha_\tau^{k-1}))E_\tau^k\Colon E_\tau^k
+\frac12\bbC(\alpha_\tau^{k-1})(E_\tau^k\Colon E_\tau^k-E_\tau^{k-1}\Colon E_\tau^{k-1})
\\&\qquad=\frac12(\alpha_\tau^k-\alpha_\tau^{k-1}))\bbC^\circ(\alpha_\tau^k,\alpha_\tau^{k-1})
E_\tau^k\Colon E_\tau^k+
\bbC(\alpha_\tau^{k-1})\Big(\frac{E_\tau^k+E_\tau^{k-1}}2\Big)\Colon(E_\tau^k-E_\tau^{k-1})\,.
\end{align*}}
\begin{align}\nonumber
&\frac12\bbC(\alpha_\tau(T))(e_{\rm el,\tau}(T)-\varepsilon_{\rm el,\tau}(T))
\Colon(e_{\rm el,\tau}(T)-\varepsilon_{\rm el,\tau}(T))
\\[-.3em]&\ \nonumber=\int_0^T\Big(\bbC(\underlinealpha_\tau)(\underline{\overline e}_{\rm el,\tau}{-}\underline{\overline\varepsilon}_{\rm el,\tau})\Colon(\DT e_{\rm el,\tau}{-}\DT\varepsilon_{\rm el,\tau})
-\frac12\DT\alpha_\tau\bbC^\circ(\overlinealpha_\tau,\underlinealpha_\tau)
(\underline{\overline e}_{\rm el,\tau}{-}\underline{\overline\varepsilon}_{\rm el,\tau})\Colon
\\[-.3em]&\qquad\qquad\qquad
\nonumber\Colon(\underline{\overline e}_{\rm el,\tau}{-}\underline{\overline\varepsilon}_{\rm el,\tau})\Big)\,\d t\ge\int_0^T\bbC(\underlinealpha_\tau)(\underline{\overline e}_{\rm el,\tau}{-}\underline{\overline\varepsilon}_{\rm el,\tau})\Colon(\DT e_{\rm el,\tau}{-}\DT\varepsilon_{\rm el,\tau})\,\d t
\end{align}
a.e.\ on $\Omega$.
When integrate over $\Omega$, this allows us to estimate:
\begin{align}\nonumber
&\int_\Omega\frac12\bbC(\alpha_\tau(T))(e_{\rm el,\tau}(T)-\varepsilon_{\rm el,\tau}(T))
\Colon(e_{\rm el,\tau}(T)-\varepsilon_{\rm el,\tau}(T))\,\d x
+\int_Q\bbD(\underlinealpha_\tau)(\DT e_{\rm el,\tau}{-}\DT\varepsilon_{\rm el,\tau})
\Colon(\DT e_{\rm el,\tau}{-}\DT\varepsilon_{\rm el,\tau})\,\d x\d t
\\[-.3em]&\nonumber\quad\le
\int_Q\Big(\bbD(\underlinealpha_\tau)(\DT e_{\rm el,\tau}{-}\DT\varepsilon_{\rm el,\tau})+
\bbC(\underlinealpha_\tau)(\underline{\overline e}_{\rm el,\tau}{-}\underline{\overline\varepsilon}_{\rm el,\tau})\Big)\Colon(\DT e_{\rm el,\tau}{-}\DT\varepsilon_{\rm el,\tau})
\,\d x\d t
\\[-.0em]&\nonumber\quad=
\int_Q(\overline f_\tau-\varrho\DT v_\tau)\Cdot(\underline{\overline v}_\tau{-}
\DT{\widetilde u}_\tau)
-\Big(\bbD(\underlinealpha_\tau)\DT\varepsilon_{\rm el,\tau}+
\bbC(\underlinealpha_\tau)\underline{\overline\varepsilon}_{\rm el,\tau}\Big)\Colon(\DT e_{\rm el,\tau}{-}\DT\varepsilon_{\rm el,\tau})
\\[-.7em]&
\qquad\
+\SYLD(\alpha_\tau)(|\DT\pi|-|\DT\pi_\tau|)
+H\underline{\overline \pi}_\tau\Colon(\DT\pi-\DT\pi_\tau)
+\kappa_1\nabla\underline{\overline\pi}_\tau{\vdots}\nabla(\DT\pi-\DT\pi_\tau)
\to 0
\label{strong-conv-plast-dam}\end{align}
From this $\DT e_{\rm el,\tau}{-}\DT\varepsilon_{\rm el,\tau}\to0$ strongly in 
$L^2(Q;\Rsym)$. Since $\DT\varepsilon_{\rm el,\tau}\to\DT e_{\rm el}$, we obtain
$\DT e_{\rm el,\tau}\to\DT e_{\rm el}$ and hence also
${\overline e}_{\rm el,\tau}\to e_{\rm el}$ strongly in 
$L^2(Q;\Rsym)$ needed for the limit passage in 
\eq{IBVP-damage-large-disc-3-plast}. The convergence in the other
terms in \eq{IBVP-damage-small-disc+plast} is then simple.

When built into the {\it phase-field 
fracture}\index{fracture!phase-field approach} model of the type 
\eq{eq6:AM-engr+}, we obtain the {\it mode-sensitive 
fracture}.\index{fracture!mode-sensitive} 
Since the fracture toughness $\Frakg_{\rm c}$ is scaled as $\mathscr{O}(1/\eps)$ 
in \eq{eq6:AM-engr+}, $\SYLD$ in \eq{plast-dam-syst-mixed-zeta}
is to be scaled as $\mathscr{O}(1/\sqrt\eps)$.
A combination of damage in its {\it phase-field 
fracture}\index{fracture!phase-field approach}
or the crack approximation 
with plasticity is referred to (an approximation of) {\it ductile cracks},
in contrast to {\it brittle cracks}\index{crack!brittle/ductile} without 
possibility of plastification on the crack tips. The idea to involve 
plastification processes into fracture mechanics is due to G.\,Irwin 
\cite{Irwi57ASSE}.

A combination of damage with the {\it perfect 
plasticity}\index{plasticity!perfect} (i.e.\ $G_\text{\sc nh}=H=0$ and 
$\kappa_1=0$) has been analysed in \cite{DaRoSt??DPPD} by using 
the strategy \eq{strong-conv-plast-dam} except that $\DDT u$ was 
approximated by the backward second time difference.


\subsection{Damage models at large strains}\label{sect-large}

In some applications, the small-strain approximation is not appropriate and
one must take into account large strains. In solid mechanics, mathematical
analysis is to be performed in the fixed reference configuration 
$\Omega\subset\R^d$, the deformation being $y:\Omega\to\R^d$. The stored
energy then depends on the deformation gradient $\nabla y$, and can be 
considered enhanced as
\begin{align}\nonumber
&\varphi_\text{\sc e}(F,\nabla F,\alpha,\nabla\alpha):=\varphi(F,\alpha)
+\mathscr{H}(\nabla F)
+\delta_{[0,1]}^{}(\alpha)+\mathscr{G}(F,\nabla\alpha)
\\&\nonumber
\qquad\COL{\text{ with }\ \ \ \mathscr{H}(G)
:=\frac14\int_{\Omega\times\Omega}\!\!
  (G(x){-}G(\wt x)):\mathbb K(x{-}\wt x)
  :(G(x){-}G(\wt x))\,\d x\d\wt x
}
\\&\label{options-gradient}
\qquad\text{ and with }\ \ \ \mathscr{G}(F,\nabla\alpha)=\frac\kappa2|\nabla\alpha|^2
\ \ \ \text{ or }\ \ \ \frac\kappa2|F^{-\top}\nabla\alpha|^2.
\end{align}
The former option in \eqref{options-gradient} is the gradient-damage theory in 
the material (reference) configuration\COL{, and is mathematically
  simpler and even a local nonsimple-material concept
  $\mathscr{H}(G)=\frac12\int_\Omega G:\mathbb K:G\,\d x$ might be used.}

The latter, \COL{mathematically more difficult} option is in the actual 
(deformed) configuration, the factor $F^{-\top}:=(F^{-1})^\top$ being the 
push-forward transformation of the vector $\nabla\alpha$ \COL{from} the
reference configuration \COL{into the actual one}. Both options are
relevant in particular situations, the latter one being mathematically more
difficult \COL{and, except \cite[Sect.\,9.5.1]{KruRou19MMCM}, has been so far
  rather devised without any rigorous proofs, cf.\ e.g.\
  \cite{Pawl06TCCH,WPMB14GELD}}.
In this latter option, the system resulted via the extended Hamilton 
variational principle \eqref{hamilton}--\eqref{Lagrangean} then reads as
\begin{subequations}\label{IBVP-damage-large+}
\begin{align}\nonumber
&\varrho\DDT y-{\rm div}\big(\pl_F^{}\varphi(\nabla y,\alpha)
+\sigma_\text{\sc k}(\nabla y,\COL{\nabla}\alpha)
-{\rm div}\,\mathfrak{H}(\nabla^2 y)\big)=f 
\\&\qquad\qquad\qquad\COL{\text{with }\
\big[\mathfrak{H}(G)\big](x)=\int_\Omega\mathbb K(x{-}\wt x)(G(x)-G(\wt x))\,\d\wt x
\ \ \text{ and}\hspace*{-2em}}
\nonumber
\\&\qquad\qquad\qquad\label{IBVP-damage-large-1+}
\text{with }\ \sigma_\text{\sc k}(F,\nabla\alpha)=
\COL{\kappa}F^{-\top}{:}(F^{-\top})'{:}(\nabla\alpha\otimes\nabla\alpha)
&&\text{ in }\ Q,
\\
&
\pl\zeta(\DT\alpha)+\pl_\alpha^{}\varphi(\nabla y,\alpha)
\ni
{\rm div}\big(\COL{\kappa}
(\nabla y)^{\COL{-1}}(\nabla y)^{-\top}\nabla\alpha\big)
\label{IBVP-damage--large-2+}
&&\text{ in }\ Q,
\\\nonumber
&
\COL{(\pl_F^{}\varphi(\nabla y,\alpha)
+\sigma_\text{\sc k}(\nabla y,\nabla\alpha))}
\vec{n}-{\rm div}_\text{\sc s}((\mathfrak{H}\nabla^2y)\COL{\cdot\vec{n}})
  =g 
\\&\qquad\qquad\qquad\qquad\qquad
\ \ \text{ and }\ \ (\mathfrak{H}\nabla^2y)\Colon(\vec{n}\otimes\vec{n})=0
&&\text{ on }\ \Sigma,
\label{IBVP-damage--large-3+}
\\&\label{IBVP-damage--large-4+}
\COL{\kappa}(\nabla y)^{\COL{-1}}(\nabla y)^{-\top}\nabla\alpha\Cdot\vec{n}=0
&&\text{ on }\ \Sigma,
\end{align}
\end{subequations}
where ${\rm div}_\text{\sc s}$ is the surface divergence. The analysis
now needs the strong convergence of $\nabla\alpha$ which now
occurs nonlinearly in the Korteweg-like stress 
$\sigma_\text{\sc k}=\sigma_\text{\sc k}(F,\nabla\alpha)$ in 
\eq{IBVP-damage-large-1+}. \COL{An important aspect is that we need to have
  a control over $(\nabla y)^{-1}$, i.e.\ $\det(\nabla y)$ should be kept
  surely away 0. As also physically desirable, this  can be ensured by
  by preventing local self-interpenetration by assuming a singulariy
  in the stored energy when $\det\,F\to0+$. More specifically,
  the potential $\varphi$ is to be qualified as 
\begin{subequations}\label{ass-growrt-phi-dam}\begin{align}
&\varphi:{\rm GL}^+(d)\times[0,1]\to\R\ \ \text{ continuously differentiable and
}
\\&\nonumber
\exists\epsilon>0\ \ \forall F\In{\rm GL}^+(d),\ \alpha\In[0,1]:\ \ 
\\[-.3em]&\qquad\varphi(F,\alpha)\ge\frac\epsilon{(\det\,F)^q}\ \ 
\text{ with }\ \ 
q>\frac{2d}{2\gamma{+}2{-}d}\ \ \text{ for some }\ \gamma>\frac d2-1\
\label{K-singular}
\intertext{while $\varphi(F,\alpha)=+\infty$ if $\det\,F\le0$,
where $\gamma$ related to the qualification of the kernel
$\mathbb K$ in \eq{options-gradient} as}
\nonumber\\[-2.4em]
&\label{K-singular+}
\exists\,\eps>0\ \forall x\In\Omega,\ F\In\R^{d\times d}:\quad
\bigg(\frac{\eps|F|^2}{|x|^{d+2\gamma}\!\!\!}\ \ -\frac1\eps\bigg)^+\!\!
\le F\Colon\mathbb K(x)\Colon F
\le\frac{|F|^2}{\eps|x|^{d+2\gamma}\!\!}\ .
\end{align}\end{subequations}
Together with the intertial term, \eq{K-singular+} grants coercivity in
$H^{2+\gamma}(\Omega;\R^d)$ which is embedded, by \eq{K-singular},
into $W^{2,p}(\Omega;\R^d)$ with $p>d$.
  
Again we use Galerkin approximation and denote the approximate solution
by $(y_k,\alpha_k)$. Testing (\ref{IBVP-damage-large+}a,b) in its
Galerkin approximation by $(\DT y_k,\DT\alpha_k)$, we obtain
 the estimates 
 \begin{subequations}\label{est-large}\begin{align}\label{est-large-1}
&\|y_k\|_{L^\infty(I;H^{2+\gamma}(\Omega;\R^d))\cap W^{1,\infty}(I;L^2(\Omega;\R^d))}\le K,\ \ \ \
\Big\|\frac1{\det(\nabla y_k)}\Big\|_{L^\infty(Q)}^{}\le K,\ \ \ \
\\&\|\alpha_k\|_{L^\infty(Q)}^{}\le K\ \ \text{ and }\ \  
\|(\nabla y_k)^{-\top}\nabla\alpha_k\|_{L^\infty(I;L^2(\Omega;\R^d))}^{}\le K\,.
\end{align}\end{subequations}
For the latter estimate in \eq{est-large-1}, we
use the result by Healey and Kr\"omer \cite{HeaKro09IWSS}, which also
excludes the Lavrentiev phenomenon.

 Based on the weak convergence $y_k\to y$ and $\alpha_k\to\alpha$ and
the Aubin-Lions compactness arguments,  we prove the 
convergence in the damage flow rule.

To prove the mentioned strong convergence $\nabla\alpha_k$, 
we use the uniform (with respect to $y$) strong monotonicity of the mapping
$\alpha\mapsto-{\rm div}\big((\nabla y)^{-1}\kappa(\nabla y)^{-\top}\nabla\alpha\big)$.
Taking $\widetilde\alpha_k$ an approximation of $\alpha$ valued in the
respective finite-dimensional spaces used for the Galerkin approximation
and converging to $\alpha$ strongly, we can test \eq{IBVP-damage--large-2+} in
its Galerkin approximation by $\alpha_k{-}\widetilde\alpha_k$ and use it in
the estimate 
\begin{align}\nonumber
&\!\!\limsup_{k\to\infty}\int_Q(\nabla y_k)^{-1}\kappa(\nabla y_k)^{-\top}
\nabla(\alpha_k{-}\widetilde\alpha_k)\Cdot\nabla(\alpha_k{-}\widetilde\alpha_k)
\,\d x\d t
\\[-.4em]&\nonumber\quad\qquad
=\lim_{k\to\infty}\int_Q\big(\pl_\alpha\varphi(\nabla y_k,\alpha_k)
+\pl\zeta(\DT\alpha_k)\big)(\widetilde\alpha_k{-}\alpha_k)
-(\nabla y_k)^{-1}\kappa(\nabla y_k)^{-\top}\nabla\widetilde\alpha_k
\Cdot\nabla(\alpha_k{-}\widetilde\alpha_k)\,\d x\d t=0
\end{align}
because $\pl_\alpha\varphi(\nabla y_k,\alpha_k)+\pl\zeta(\DT\alpha_k)$ is
bounded in $L^2(Q)$ while $\widetilde\alpha_k{-}\alpha_k\to0$ strongly in
$L^2(Q)$ by the Aubin-Lions compactness theorem and because 
$(\nabla y_k)^{-1}\kappa(\nabla y_k)^{-\top}\nabla\widetilde\alpha_k$ converges
strongly in $L^2(Q;\R^d)$ while $\nabla(\alpha_k{-}\widetilde\alpha_k)\to0$
weakly in $L^2(Q;\R^d)$. As $(\nabla y_k)^{-1}\kappa(\nabla y_k)^{-\top}$ is
uniformly positive definite, we thus obtain that
$\nabla(\alpha_k{-}\widetilde\alpha_k)\to0$ strongly in $L^2(Q;\R^d)$,
and thus $\nabla\alpha_k\to\nabla\alpha$ strongly in $L^2(Q;\R^d)$.
Then we have the convergence in the Korteweg-like stress
$\sigma_\text{\sc k}(\nabla y_k,\nabla\alpha_k)
\to
\sigma_\text{\sc k}(\nabla y,\nabla\alpha)$
even strongly in $L^p(I;L^1(\Omega;\R^{d\times d}))$ for any $1\le p<+\infty$.
 The limit passage in the force equilibrium towards 
\eq{IBVP-damage-large-1+} formulated weakly is then straightforward. 

}

{\small
  \bigskip\bigskip \noindent{\it Acknowledgement}
Special thanks are to Roman Vodi\v{c}ka for providing sample snapshots 
from numerical simulations presented in Fig.\:\ref{fig-simul}. 
Also discussion with Vladislav Manti\v c about finite fracture mechanics
and coupled criterion has been extremely useful, as well as discussions 
with
Martin Kru\v z\'\i k about the 
\COL{actual} \NEW{gradient of damage} at large strains.
\COL{Careful reading of the 
manuscript and many comments of Elisa Davoli and Roman Vodi\v cka, as well as
of an anonymous referee are also appreciated very much.} A partial support from
the Czech Science Foundation projects
17-04301S and \COL{19-04956S}, 
the institutional support RVO:\,61388998 (\v CR), and also 
by the Austrian-Czech project 16-34894L (FWF/CSF) are acknowledged, too.
}

\COMMENT{CHECK ALSO 






Gradient damage modeling of brittle fracture in an explicit
dynamics context
Tianyi Li and Jean-Jacques Marigo and Daniel Guilbaudand Serguei Potapov
INTERNATIONAL JOURNAL FOR NUMERICAL METHODS IN ENGINEERING
Int. J. Numer. Meth. Engng 2016; 108:1381-1405

Variational Approach to Dynamic Brittle Fracture
via Gradient Damage Models
Tianyi Li and Jean-Jacques Marigo and Daniel Guilbaud and Serguei Potapov
Applied Mechanics and Materials Submitted: 2015 Vol. 784, pp 334-341
NOTE: staggered scheme + Newmark 

A length scale insensitive phase-field damage model for brittle fracture
Jian-Ying Wu and Vinh Phu Nguyen
Journal of the Mechanics and Physics of Solids 
119, 2018, Pages 20--42


Wu, J. Y., 2018a. A geometrically regularized gradient-damage model with 
energetic equivalence. Comput. Methods Appl. Mech. Engrg. 328,
612--637.

Wu, J. Y., 2018b. Numerical implementation of non-standard phase-field damage 
models. Computer Methods in Applied Mechanics and Engineering in press.

ALSO \cite{BVSH12PFDD} AND \cite{HofMie13PFMD} = M Hofacker, C Miehe:
A phase field model of dynamic fracture: Robust field updates for the analysis
of complex crack patterns. 
Int J Numer Methods Eng 93 (2013):276--301

}

\end{sloppypar}
\end{document}